\documentclass[12pt]{article2}
\usepackage{amsmath}
\usepackage{amssymb}
\usepackage{amsthm}
\usepackage{graphicx}
\numberwithin{equation}{section}
\newtheorem{theorem}{Theorem}[section] 
\newtheorem{proposition}[theorem]{Proposition}
\newtheorem{lemma}[theorem]{Lemma}

\newtheorem{Definition}[theorem]{Definition}
\newenvironment{definition}{\begin{Definition}\rm}{\end{Definition}}
\newtheorem{Remark}[theorem]{Remark}
\newenvironment{remark}{\begin{Remark}\rm}{\end{Remark}}
\newtheorem{Example}[theorem]{Example}

\newcommand\gog{\mathfrak{g}}
\newcommand\CC{\mathbb{C}}
\newcommand\ZZ{\mathbb{Z}}
\newcommand\Zplus{\mathbb{Z}_{\geq 0}}

\newcommand\ga{\gamma}
\newcommand\de{\delta}
\newcommand\ep{\varepsilon}
\newcommand\la{\lambda}

\newcommand\De{\Delta}
\newcommand\Om{\Omega}
\newcommand\FSA{{\cal A}}
\newcommand\FSC{{\cal C}}
\newcommand\FSF{{\cal F}}
\newcommand\FSM{{\cal M}}
\newcommand\FSR{{\cal R}}
\newcommand\FSU{{\cal U}}
\newcommand\cd{\cdot}
\newcommand\half{\frac12}
\newcommand\thalf{{\tfrac12}}
\newcommand\quart{\frac14}
\newcommand\tquart{\tfrac14}
\newcommand\Uq{\FSU_q}
\newcommand\id{{\rm id}}
\newcommand\ten{\otimes}
\newcommand\iy{\infty}
\newcommand\wt{{\rm wt}}
\newcommand{\qhyp}[5]{\,\mbox{}_{#1}\phi_{#2}\left(
  \genfrac{}{}{0pt}{}{#3}{#4};#5\right)}
\newcommand{\hyp}[5]{\,\mbox{}_{#1}F_{#2}\left(
  \genfrac{}{}{0pt}{}{#3}{#4};#5\right)}
\newcommand\qbinom[3]{\genfrac[]{0pt}0#1#2_#3}
\newcommand\Proof{\smallskip\noindent{\bf Proof}\quad}
\newcommand\lan{\langle}
\newcommand\ran{\rangle}
\newcommand\LHS{left-hand side}
\newcommand\RHS{right-hand side}
\newcommand\qpoch[2]{([#1]_q)_{#2}}
\newcommand\qfac[1]{[#1]_q!}
\newcommand\qthreej[6]{\left[ \begin{array}{lcr}#1&#2&#3\\#4&#5&#6\\
 \end{array}\right]_{q}}
\newcommand\qsixj[6]{\left\{ \begin{array}{lcr}#1&#2&#3\\#4&#5&#6\\
 \end{array}\right\}_{q}}
\newcommand\wtilde{\widetilde}

\begin{document}

\title{Fusion and exchange matrices for quantized
sl(2) and associated q-special functions}
\author{T.~H. Koornwinder and N. Touhami}

\date{February 19, 2003}

\maketitle
\section{Introduction}
The aim of this paper is to evaluate in terms of $q$-special functions
the objects (intertwining map, fusion matrix, exchange matrix)
related to the quantum
dynamical Yang-Baxter equation (QDYBE) for infinite dimensional 
representations (Verma modules)
of the quantized universal enveloping algebra
$\Uq(\gog)$ in the case $\gog=sl(2,\CC)$.
This study is done in the framework of the  exchange construction,
initiated by Etingof and Varchenko, in order to find solutions to 
the QDYBE for $\Uq(\gog)$ with $\gog$ a complex semisimple Lie algebra
(see \cite[\S2]{EV}, \cite[\S2]{ES}, \cite[\S3]{E}).
As a result, the familiar
interpretation
of $q$-Hahn and $q$-Racah polynomials as $q$-Clebsch-Gordan and $q$-Racah
coefficients, respectively, for finite dimensional irreducible
representations of
the quantum group $SU_q(2)$ is extended to a much larger range of parameters,
while moreover these interpretations are obtained in an unusual and
interesting way, following the definitions of intertwining map and
exchange matrix. Furthermore, we reprove, by using these
explicit expressions, some properties related to the objects under study
in the special case $\gog=sl(2)$, which were earlier proved
in \cite{EV} and \cite{ES} in a more
abstract way in the case of more general $\gog$.

The paper is organized as follows.
In Section~2 a brief review of $q$-special functions and
of the quantized universal algebra $\Uq(sl(2))$
and its representations is given.
Section~3
deals with the intertwining operator 
$\Phi_{q,\la}^{v}: M_{q,\la} \longrightarrow M_{q,\mu}\ten V$.
Here $M_{q,\la}$, $M_{q,\mu}$ are Verma modules for $\Uq(sl(2))$
with highest weight vectors $x_\la,x_\mu$ ($\la,\mu\in\CC$),
$V$ is a $\Uq (sl(2))$-module, not necessarily of finite dimension and
soon assumed to be a Verma module, $\Phi_{q,\la}^{v}$ is
$\Uq(sl(2))$-intertwining,
and $\Phi_{q,\la}^{v} (x_{\la})$ is supposed to have ``highest''
term $x_\mu\ten v$.
On various places in the paper we do explicit computations first in
the generic infinite dimensional case and next make transition
(by continuity)
to the finite dimensional case in order to connect the results with
known explicit expressions in the finite dimenisonal case.
These limit transitions need careful justification.

The matrix elements of the intertwining operator with respect to
the standard bases of 
Verma modules generalize the $q$-Clebsch-Gordan coefficients for
the finite dimensional case and
can still be expressed 
in terms of $q$-Hahn polynomials
(see for example \cite{BiLoh} for the finite dimensional case).
This result is used in Section~4 to obtain in a new way the
(otherwise known) orthogonality 
relations of $q$-Hahn polynomials by use of the $q$-analogue of the
Shapovalov form.
In Section~5 we find the explicit matrix elements of the 
fusion matrix (and of its inverse)
$J_{W,V}: W\ten V \longrightarrow V\ten W$,
which is defined in terms of the intertwining operator by
the identity
$(\Phi_{q,\la - wt(v)}^{w}\ten \id)\circ\Phi_{q,\la}^{v}
=\Phi_{q,\la}^{J_{W,V}(\la)(w \ten v)}$.
This leads to the expression of the universal fusion matrix in Section 6:
a generalized
element of $\Uq(sl(2)) \ten \Uq(sl(2))$ which is shown
to satisfy a shifted 2-cocycle condition (results
earlier observed in \cite{BabBB}).
We know also from
\cite{BabBB} that one can associate to this element a   
generalized element of $\Uq(sl(2))$ which is called the shifted boundary.
We derive the explicit
expression for its inverse independently in Section~7 and we  point out
how this indeed can be seen as the inverse of the shifted boundary in
\cite{BabBB}.
We also derive the result, first observed by Rosengren
\cite{R2000b}, that the shifted
boundary acts as a generalized conjugation in $\Uq(sl(2))$
which sends a standard $q$-analogue of a Cartan subalgebra to
e non-standard $q$-analogue of a Cartan subalgebra.
Next, in Section 8, we derive the 
ABRR equation for the universal fusion matrix (earlier
observed in \cite{ArBRR} and \cite{ES} for more general $\gog$). 

In Section~9 we compute the matrix elements of the exchange matrix
$R_{q,\ga,\de}(\la)$ sending $M_{q,\ga} \ten M_{q,\de}$
into itself. These turn out to be
$q$-Racah polynomials.
In Section~10 we show that these matrix elements in the finite dimensional
case are essentially $q$-Racah coefficients ($q$-$6j$ symbols).
The exchange matrix is known to satisfy the QDYBE.
In Section~11 we prove in more detail the observation made in
\cite[\S 8]{EV} that
the QDYBE in the finite dimensional
case is equivalent to a known
identity (see \cite{KiR} and 
\cite{N1989}) satisfied by $q$-Racah coefficients ($q$-$6j$ symbols), which is
called the Yang-Baxter equation for the 
interaction round a face (IRF) model (see \cite{Bax}).

Note that in the classical limit $q\to1$ of the objects we are dealing with 
in this paper,
the role of
$\Uq(sl(2))$ is taken over by the Lie algebra $sl(2)$, see
\cite [\S 2.1]{ES} and \cite{KooT}.
The different expressions corresponding to the classical limit are 
obtained by taking the $q\rightarrow 1$ limit and replacing 
the $q$-hypergeometric 
functions ${}_r\Phi_s$ by their classical analogues ${}_rF_s$.
This limit can usually be taken in a straightforward way, and the
$q$-case is only computationally more difficult than the $q=1$ case.
However, in Section 9
the derivation of the exchange matrix in the $q$-case requires
one more summation than in the $q=1$ case. This is because the definition
of the exchange matrix additionally involves the $R$-matrix in the $q$-case.
This also complicates things a little bit in Section 10.
The only place in this paper where the limit transition $q\rightarrow 1$
fails, is for
the shifted boundary \eqref{eq:SB2}. This object does not seem to exist
for $q=1$.
In the case of the exchange matrix in the 
finite-dimensional case \eqref{eq:ER22} the classical limit can be 
related to Racah coefficients and $6j$ symbols, see the definitions  and  
classical analogues  to \eqref{eq:ER9}--\eqref{eq:ER11} in  
\cite[\S 3.18]{BiLou}. 
Then, similarly as in Section 11, one can show also for $q=1$ that the QDYBE
yields a known (see \cite{N1988}) 
identity for sums of products of three $6j$ symbols.
\\[\medskipamount]\goodbreak\noindent
{\sl Note}\quad
Some of the material of this paper was earlier presented by the first author
\cite{Koo2000b} at the Advanced Study Institute 
Special Functions 2000, Arizona State University.
\\[\medskipamount]\goodbreak\noindent
{\sl Conventions}\quad
Throughout this paper we assume that $0<q<1$.\\
We use the notations 
$m\wedge n:=\min(m,n)$ and $m\vee n:=\max(m,n)$.
\section{Preliminaries} \label{sec:PR}
\subsection{q-Hypergeometric functions}
\label{subsec:PR1}
Standard references are
Gasper and Rahman \cite{GR} and (for special orthogonal polynomials)
Koekoek and Swarttouw \cite{KoeS}.
We recall the definition and notation for the $q$-Pochhammer symbol and
the $q$-binomial coefficient which is nowadays standard in the theory
of $q$-special functions
(see \cite{GR}):
\begin{eqnarray}
&&\qquad\quad\;\;\;(a;q)_{k}:= (1-a) (1-aq) \ldots (1-aq^{k-1})\qquad
(k\in\Zplus),
\label{eq:PR1}\\[\medskipamount]
&&\qquad\quad\;\;(a;q)_\iy:=\textstyle\prod_{j=0}^\iy\,(1-aq^j)\nonumber\\
&&(a_1, \ldots, a_n;q)_k:= (a_1;q)_k \ldots (a_n;q)_k\,,\qquad
{\qbinom nkq}:=\frac{(q;q)_n}{(q;q)_k\,(q;q)_{n-k}}\,.\nonumber
\end{eqnarray}
The following, slightly different definition and notation for the $q$-number,
the $q$-factorial and the $q$-Pochhammer symbol
is very convenient for computations in quantum groups:
\begin{equation}
 [a]_q :=\frac{q^{a/2}-q^{-a/2}}{q^{1/2}-q^{-1/2}}\,,\quad
	\qfac k := \textstyle\prod_{j=1}^k\, [j]_q,\quad
	\qpoch a k := \prod_{j=0}^{k-1}\,[a+j]_q
	\qquad(k\in\Zplus).
\label{eq:PR24}
\end{equation}
The symbols $\qfac k$ and $\qpoch a k$ can be expressed
in terms of the standard notation \eqref{eq:PR1}
as follows:
\begin{equation}
\qfac k =q^{-\quart k(k-1)} \frac{(q;q)_{k}}{(1-q)^{k}}\,,\qquad
\qpoch a k=q^{-\half k(a-1)} q^{-\quart k(k-1)} 
	{\frac{(q^{a};q)_{k}}{(1-q)^{k}}}\,.
\label{eq:PR2}
\end{equation}
The following consequences of these identities are also useful:
\begin{equation}
\frac{\qpoch a k}{\qpoch b k}=q^{-\half k(a-b)}\frac{(q^a;q)_k}{(q^b;q)_k}\,,
\qquad
\frac{\qfac n}{\qfac k\,\qfac{n-k}}=q^{-\half k(n-k)}\qbinom nkq.
\label{eq:PR25}
\end{equation}

The general {\em $q$-hypergeometric series} is given by
\begin{equation}
\qhyp rs{a_1,\ldots,a_r}{b_1,\ldots,b_s}{q,z}
:= \sum_{k=0}^\iy \frac{(a_1,\ldots,a_r;q)_k}
{(b_1,\ldots,b_s;q)_k\,(q;q)_k}\,\bigl((-1)^k
q^{\half k(k-1)}\bigr)^{s-r+1} z^k,
\label{eq:PR23}
\end{equation}
which is an infinite series if none of $a_1,\ldots,a_r$ is equal to
$q^{-n}$ for some $n\in\Zplus$ and is terminating otherwise.
In the infinite series case we require for convergence
that $r<s+1$ or that $r=s+1$ and $|z|<1$.
The coefficients $b_j$ are not allowed to be equal to $q^{-m}$ for some
$m\in\Zplus$ unless for some $i$ we have $a_i=q^{-n}$ with $n\in\Zplus$
and $n\le m$.
If in \eqref{eq:PR23} we replace $a_i,b_j$ by
$q^{a_i},q^{b_j}$, respectively, and if we let
$q\uparrow 1$ then the limit equals the
{\em hypergeometric series}
\begin{equation}
\hyp rs{a_1,\ldots, a_r}{b_1,\ldots, b_s}z:=
\sum_{k=0}^\iy\frac{(a_1)_k\ldots(a_r)_k}
{(b_1)_k\ldots(b_s)_k\,k!}\,z^k,
\label{eq:PR26}
\end{equation}
where $(a)_k:=a(a+1)\ldots(a+k-1)$ is the Pochhammer symbol.

We will need some summation formulas, namely
the {\em $q$-Chu-Vandermonde sum} \cite[(II.6)]{GR}
\begin{equation}
	\qhyp21{q^{-n},a}c{q,q} = a^n\, \frac{(c/a;q)_n}
	{(c;q)_n}\qquad(n\in\Zplus),
\label{eq:PR27}
\end{equation}
a variant \cite[(II.7)]{GR} of the
$q$-Chu-Vandermonde sum (obtained by reverting the order of summation in
\eqref{eq:PR27})
\begin{equation}
	\qhyp21{q^{-n},a}c{q,\frac{q^n c}a} = 
	 \frac{(c/a;q)_n}{(c;q)_n}\qquad(n\in\Zplus),
\label{eq:PR28}
\end{equation}
the limit of \eqref{eq:PR28} for $c\rightarrow \iy$ given by
\begin{equation}
\qhyp20{q^{-n},a}{-}{q,\frac{q^n}a}=a^{-n}\qquad(n\in\Zplus),
\label{eq:PR29}
\end{equation}
and the limit of \eqref{eq:PR28} for $n\rightarrow \iy$
(see \cite[(II.5)]{GR}):
\begin{equation}
	\qhyp11{a}{c}{q,\frac ca} =  \frac{(c/a;q)_{\iy}}{(c;q)_{\iy}}\,.
\label{eq:PR35}
\end{equation}

We will also need  some transformation formulas, namely
a  ${}_3\phi_2$ transformation 
\cite[(10.10.5)]{AnAR}
\begin{equation}
	\qhyp32{q^{-n},q^a,q^b}{q^d,q^e}{q,q}
	=
	q^{an}\frac{(q^{e-a};q)_n}{(q^e;q)_n} 
	\qhyp32{q^{-n},q^a,q^{d-b}}{q^d,q^{a-e-n+1}}{q,q}\qquad
(n\in\Zplus),
\label{eq:PR30}
\end{equation}
and {\em Sear's transformation} \cite[(III.16)]{GR}
of a terminating balanced ${}_4\phi_3$ series:
\begin{eqnarray}
	\qhyp43{q^{-n},q^a,q^b,q^c}{q^d,q^e,q^f}{q,q}
	= \frac{(q^a, q^{e+f-a-b},q^{e+f-a-c};q)_n}
	{(q^e, q^f,q^{e+f-a-b-c};q)_n }
	\qhyp43{q^{-n},q^{e-a},q^{f-a},q^{e+f-a-b-c}}
	{q^{e+f-a-b},q^{e+f-a-c},q^{1-n-a}}{q,q}
\nonumber&\\[\medskipamount]
(n\in\Zplus,\quad a+b+c-n+1 = d+e+f).\qquad\qquad
\label{eq:PR31}
\end{eqnarray}

We will meet some $q$-hypergeometric orthogonal polynomials,
namely {\em $q$-Hahn polynomials}
\begin{equation}
Q_n(x;a,b,N;q):=
\qhyp32{q^{-n},q^{n+1}ab,q^{-x}}{qa,q^{-N}}{q,q}\qquad
(N\in\Zplus,\;n\in\{0,1,\ldots,N\})
\label{eq:PR40}
\end{equation}
(see \cite[\S 3.6]{KoeS}),
and {\em $q$-Racah polynomials}, which are defined
in \eqref{eq:EM7}
in terms of certain terminating balanced ${}_4\phi_3$ series,
(see also \cite[\S 3.2]{KoeS}).

In section \ref{sec:SB} we will use
the $q$-analogues of the exponential function
(see \cite{GR}) given by
\begin{equation}
	e_q(z) := \sum_{n=0}^{\iy} \frac{z^n}{(q;q)_n}=\frac1{(z;q)_\iy}\,, 
	\qquad E_q(z):= \sum_{n=0}^{\iy} \frac{q^{\half n(n-1)}z^n}{(q;q)_n}=
(-z;q)_\iy\,,
\label{eq:PR36}
\end{equation}
which are related (for $|z|<1$ or as
formal power series) by
\begin{equation}
	e_q(z) E_q(-z) = 1,
\label{eq:PR37}
\end{equation}
and we will also use there a formal power series in noncommuting variables
\begin{equation}
 	\FSA_q(x,y) = \sum_{k=0}^{\iy} \frac{1}{(q;q)_k}\, y^k\, \frac{1}{(x;q)_k}
\qquad(xy = q^2 yx),
\label{eq:PR38}
\end{equation}
which is invertible
(see \cite[Lemma 3.4]{R2000a}) with the following inverse:
\begin{equation}
\FSA^{-1}_q(x,y) =
\sum_{k=0}^{\iy} \frac{(-1)^k q^{\half k(k-1)}}{(q;q)_k}\,
\frac{1}{(q^{-k-1}x;q)_k}\,y^k.
\label{eq:PR39}
\end{equation}
\subsection{Preliminaries on $\Uq(sl(2,\CC))$}
We use Chari and Pressley \cite{CP} as a standard reference for quantum
groups.
The quantized universal enveloping algebra $\Uq:=\Uq(sl(2,\CC))$ 
is the algebra generated by the elements $e$, $f$, $q^{h/4}$,
$q^{-h/4}:=(q^{h/4})^{-1}$ satisfying relations
\begin{equation}
 q^{h/4} e = q^{1/2}\, e q^{h/4},\quad
 q^{h/4} f= q^{-1/2}\, f q^{h/4},\quad
 [e,f] = [h]_q\,.
\label{eq:PR3}
\end{equation}
It follows by induction from \eqref{eq:PR3} that
\begin{eqnarray}
	q^{\xi h} e^j = e^j q^{\xi (h+2j)},\qquad
	q^{\xi h} f^j = f^j q^{\xi (h-2j)},\qquad
	ef^{n} = f^{n} e + [n]_{q} f^{n-1}[h-n+1]_{q}\,,
\label{eq:PR18}
\end{eqnarray}
and (see for instance \cite[(1.3.1)]{DeCK})
\begin{equation}
e^mf^n=\sum_{k=0}^{\min(m,n)}
\frac{\qfac m\, \qfac n}{\qfac k\, \qfac{m-k}\, \qfac{n-k}}\,
f^{n-k}e^{m-k}\,\qpoch{h+m-n-k+1}k\,.
\label{eq:PR34}
\end{equation}
It is convenient to allow $h$ as a formal element of $\Uq$.
Now \eqref{eq:PR3} can equivalently be written as:
\begin{equation}
[h,e] = 2e,\quad
[h,f] = -2f,\quad
[e,f] = [h]_{q}.
\label{eq:PR33}
\end{equation}

$\Uq$ is endowed with a structure of 
quasi-triangular Hopf algebra (see for example \cite{CP}).
The coproduct $\De$, the counit $\ep$ 
and the antipode $S$ are given by:
\begin{equation}	
\De(h) = h \ten 1 + 1 \ten h,\quad
\De(e)  = e\ten q^{h/4} + q^{-h/4} \ten e,\quad
\De(f)  = f\ten q^{h/4} + q^{-h/4} \ten f,
\label{eq:PR4}
\end{equation}
\[
	\ep(h) = \ep (e) = \ep (f) = 0,
\]
\[
	S(h) = -h,\quad
	S(e) = - q^{1/2}e,\quad
	S(f) =- q^{-1/2}f.
\]
Note that, by the $q$-binomial theorem,
\begin{eqnarray}
	\De(e^{n}) = \sum_{k=0}^{n}
		{\frac{[n]_{q}!}{[k]_{q}!\,[n-k]_{q}!}}e^{k}q^{-(n-k)h/4}
		\ten e^{n-k} q^{k h/4},\label{{eq:PR5}}\\
	\De(f^{n}) = \sum_{k=0}^{n} 
		{\frac{[n]_{q}!}{[k]_{q}!\,[n-k]_{q}!}}f^{k}q^{-(n-k)h/4}
		\ten f^{n-k} q^{k h/4}.\label{{eq:PR6}}
\end{eqnarray}

The universal $\FSR$-matrix is given by the following expression,
see Drinfel'd \cite{Dr}:
\begin{equation}
	\FSR= q^{\quart(h \ten h)}\sum_{j=0}^\iy \frac{(1-q^{-1})^j 
	q^{-\quart j(j-1)}}{[j]_{q}!}
	(q^{\quart jh}e^{j}\ten q^{-\quart jh} f^{j}).
	\label{eq:PR7}
\end{equation}
It satisfies:
\begin{equation}
	(\id \ten \De) \FSR 
	=
	\FSR_{13} \FSR_{12},\qquad
	(\De \ten \id ) \FSR 
	=
	\FSR_{13} \FSR_{23},\qquad
	(S \ten \id ) \FSR 
	=
	\FSR^{-1}.
\label{eq:PR22}
\end{equation}
The Casimir element $\Om$ of $\Uq$ (a central element) is given by
\begin{equation}
  	\Om := fe + \left[\thalf h \right]_{q}
	\left[\thalf h+1\right]_{q}.
	\label{eq:PR8} 
\end{equation}

For a $\Uq$-module $V$ and for $\la\in\CC$ let $V[\la]$ denote
the {\em weight space
of weight $\la$} in $V$, i.e.,
\begin{equation}
	V[\la]:= \{ v \in V \mid h\cd v = \la v\}.
\label{eq:PR9}
\end{equation}
The nonzero elements of $V[\la]$ are called {\em weight vectors of
weight $\la$}.
We say that a weight $\la$ {\em occurs in $V$}
if $V[\la]\ne\{0\}$.
A weight vector $v\in V$ is called a {\em highest weight vector} if
$e\cd v=0$.
\\[\medskipamount]\goodbreak\noindent
{\sl Convention}\quad
Each $\Uq$-module to be considered will be spanned by its weight vectors
and the real parts of its occurring weights will be bounded
from above.\\[\medskipamount]\indent
Let $M_{q,\la}$ be
the {\em Verma module} for $\Uq$ with {\em highest weight vector} $x_\la$ of
{\em highest weight} $\la\in\CC$, i.e., $h\cd x_\la=\la\,x_\la$,
$e\cd x_\la=0$.
A basis of weight vectors for $M_{q,\la}$ is given by the elements 
$x^{\la}_{\la-2k}:=f^{k}\cd x_{\la}$ 
($k\in\Zplus$).
The action of the generators of $\Uq$ on this basis is given by:
\begin{eqnarray}
	h\cd x^{\la}_{\la - 2k} 
	&=& 
	(\la - 2k)\, x^{\la}_{\la -2k}, 
	\label{eq:PR10}\\
	e\cd x^{\la}_{\la - 2k}
	&=& 
	[k]_{q} [\la - k+1]_{q}\,
	x^{\la}_{\la -2k+2},
	\label{eq:PR11}\\
	f\cd x^{\la}_{\la - 2k} 
	&=& 
	 x^{\la}_{\la -2k-2}.
	\label{eq:PR12}
\end{eqnarray}
Hence,
\begin{eqnarray}
	e^j\cd x^{\la}_{\la - 2k}
	 &=& 
	(-1)^j ([-k]_{q})_j ([\la - k+1]_{q})_j\,
	x^{\la}_{\la -2k+2j},
\label{eq:PR13}\\
	f^j\cd x^{\la}_{\la - 2k} 
	&=& 
	x^{\la}_{\la -2k-2j},
	\label{eq:PR14}
\end{eqnarray}

$M_{q,\la}$ is an irreducible $\Uq$-module iff $\la\notin\Zplus$.
In the case that $\la\in\Zplus$ we see that $M_{q,-\la-2}$ (here spanned by
$x_{\la-2k}^\la$, $k>\la$) is an
irreducible submodule and that the quotient module
$M'_{q,\la}:=M_{q,\la}/M_{q,-\la-2}$ is
a finite dimensional irreducible $\Uq$-module. We can realize $M'_{q,\la}$
with basis
\begin{equation}
x^\la_{\la-2k}:=f^k\cd x_\la\qquad (k=0,1,\ldots,\la)
\label{eq:PR16}
\end{equation}
and with $\Uq$-action \eqref{eq:PR10}--\eqref{eq:PR12} except that
\begin{equation}
f\cd x_{-\la}^\la=0.
\label{eq:PR17}
\end{equation}
\begin{remark} \label{th:PR20}
The following notations for a finite dimensional irreducible
$\Uq$-module are often used in literature
(see for instance \cite{KiR}).
If $j\in\thalf\Zplus$ then write $V^j$ for
$M'(q,2j)$ and use as a standard basis for $V^j$ the vectors
\begin{equation}
e_m^j:=\frac1{\sqrt{\qpoch{j+m+1}{j-m}\qfac{j-m}}}\,x_{2m}^{2j}\qquad
(m\in\{-j,-j+1,\ldots,j\}).
\label{eq:PR21}
\end{equation}
The symmetric bilinear form $\lan\;\ran$ on $V^j$ for which
the $e_m^j$ are orthonormal (i.e., $\lan e_m^j,e_n^j\ran=\de_{m,n}$),
will satisfy \eqref{eq:OS1}, by which $\lan\;\ran$ will be a
Shapovalov form on $V^j$, to be discussed in Section \ref{sec:OS}
\end{remark}

The following universal property of Verma modules $M_{q,\la}$ will be useful.
It follows immediately from \eqref{eq:PR18} and
\eqref{eq:PR10}--\eqref{eq:PR12}.
\begin{lemma} \label{th:PR19}
Let $\la\in\CC$. Let $V$ be a $\Uq$-module with a highest weight vector
$v\in V[\la]$. Then there is a unique $\Uq$-intertwining operator
$\Phi: M_{q,\la}\to V$ such that $\Phi(x_\la)=v$. It is given by
$\Phi(f^k\cd x_\la):=f^k\cd v$. If $\la\in\Zplus$ then $\Phi$ is also
well-defined on the quotient $M'_{q,\la}$ iff $f^{\la+1}\cd v=0$.
\end{lemma}

By \eqref{eq:PR8} the Casimir operator acting on $M_{q,\la}$ is
a constant multiple of the identity:
\begin{equation}
	\Om= \left[{\thalf\la}\right]_{q}
	\left[\thalf\la +1 \right]_{q}.
	\label{eq:PR15}	
\end{equation}

The action of $\Uq$  on the tensor product of two 
$\Uq$-modules $V$ and $W$ is given by the comultiplication:
\[
x\cd(v\ten w):=\De(x)\cd(v\ten w) \qquad(x\in\Uq).
\]
Let $P: v\ten w\mapsto w\ten v: V\ten W \to W\ten V$
be the flip operator. Then
$P\circ\FSR: v\ten w\mapsto P(\FSR\cd(v\ten w)):
V\ten W \to W\ten V$
is an intertwining operator of $\Uq$-modules.
\section{The intertwining map} \label{sec:IM}
The statement and proof of \cite[Proposition 2.1]{ES}
(existence and uniqueness of the intertwining map)
can be adapted
to the quantum case and to the case that the module $V$ is not
necessarily of finite dimension. We start with a lemma
(we will only deal with $\gog=sl(2)$).
\begin{lemma} \label{th:IM2}
Let $V$ be a $\Uq$-module, let $\la,\mu\in\CC$,
and let $0\ne v\in V[\la-\mu]$.
In the case that $\mu\in\Zplus$ assume moreover that $e^{\mu+1}\cd v=0$.
Then there is a unique highest weight vector of weight $\la$ in
$M_\mu\ten V$ of the form
\begin{equation}
\sum_{k=0}^\iy f^k\cd x_\mu\ten v_k\quad
\mbox{\rm such that $v_0=v$ and such that $v_k=0$ if $k>\mu$ and
$\mu\in\Zplus$.}
\label{eq:IM3}
\end{equation}
Furthermore, there is $l\in\Zplus$ such that $v_k\ne0$ for $k=0,\ldots,l$
and $v_k=0$ for $k>l$. Then $v_l$ is a highest weight vector. Finally,
\begin{equation}
v_k=\frac{q^{-\quart k(\la + 2)}}
{\qfac k\,\qpoch{-\mu}k}\, e^{k}\cd v\qquad
\mbox{\rm($k\in\Zplus$ or $k=0,\ldots,\mu$ if $\mu\in\Zplus$),}
\label{eq:IM7}
\end{equation}
\end{lemma}
\Proof
If $v_k\ne0$ then $v_k$ must have weight $\la-\mu+2k$. By the convention
about $\Uq$-modules in Section \ref{sec:PR}, $v_k=0$ for $k$
sufficiently large.
Hence the condition $e\cd w=0$ holds iff
\begin{equation}
[k+1]_q[\mu-k]_q v_{k+1}=-q^{-\quart(\la+2)}e\cd v_k\quad(k\in\Zplus).
\label{eq:IM8}
\end{equation}
This proves existence and uniqueness of the highest weight vector $w$
of the form \eqref{eq:IM3}. (Note that in the case $\mu\in\Zplus$
the case $k=\mu$ of \eqref{eq:IM8} just says that $e\cd v_\mu=0$, and that
this follows from the assumption $e^{\mu+1}\cd v=0$ together with the cases
$k<\mu$ of \eqref{eq:IM8}.).

For the other statements of the lemma let
$l\in\Zplus$ be maximal such that $v_k=0$ for $k=0,\ldots,l$.
First assume that $\mu\notin\Zplus$.
Then it follows by iteration of \eqref{eq:IM8} that \eqref{eq:IM7}
holds for $k\in\Zplus$ and that $e\cd v_l=0$ and $v_k=0$ for $k>l$.
Next assume that $\mu\in\Zplus$.
Then it follows by iteration of \eqref{eq:IM8} that \eqref{eq:IM7}
holds for $k=0,\ldots,\mu$. Then $e\cd v_l=0$ by \eqref{eq:IM8} if $l<\mu$,
and $e\cd v_l=0$ by assumption if $l=\mu$. If $l<\mu$ then we see also
that $v_k=0$ for $k=l+1,\ldots,\mu$.
$\square$
\begin{definition} \label{th:IM5}
Let $V,\la,\mu,v$ and further assumptions be as in Lemma \ref{th:IM2}.
The {\em intertwining map} $\Phi_{q,\la}^v$ is defined as the
unique $\Uq$-intertwining operator
$\Phi_{q,\la}^v: M_{q,\la}\to M_{q,\mu}\ten V$ such that
$\Phi_{q,\la}^v(x_\la)$ is of the form \eqref{eq:IM3}.
\end{definition}
The existence and uniqueness of $\Phi_{q,\la}^v$ follow from
Lemma \ref{th:IM2} together with Lemma \ref{th:PR19}.
Then we obtain by \eqref{eq:IM7} that
\begin{equation}
\Phi_{q,\la}^{v}(x_{\la})= 
\sum_{k=0}^{\iy\;{\rm or}\;\mu} \frac{q^{-\quart k(\la + 2)}}
{\qfac k\,\qpoch{-\mu}k}\, f^{k}\cd x_{\mu} \ten e^{k}\cd v 
\label{eq:IM4} 
\end{equation}
\begin{remark} \label{th:IM9}
Consider Definition \ref{th:IM5} in the case that $\mu\in\Zplus$.
Then the sum in \eqref{eq:IM4} has upper limit $\mu$.
Suppose that moreover $V=M_{q,\ga}$ and that $\la,\ga\in\Zplus$.
Then, since $v\ne 0$ and $\wt(v)=\la-\mu$, we have $\la-\mu\le\ga$
and $\la-\mu-\ga$ is even.
Also, the condition $e^{\mu+1}\cd v=0$ in Lemma \ref{th:IM2} is certainly
satisfied if $\la+\mu\ge\ga$.

Since the canonical maps $M_{q,\mu}\to M'_{q,\mu}$ and
$M_{q,\mu}\to M'_{q,\mu}$ are $\Uq$-intertwining,
we can consider \eqref{eq:IM4} for the intertwining map
$\Phi_{q,\la}^v: M_{q,\la}\to M'_{q,\mu}\ten M'_{q,\ga}$.
In order to have the canonical projection of $v$ nonzero, we require that
$\la-\mu\ge-\ga$.
By Lemma \ref{th:PR19}, this last map induces an intertwining map
$\Phi_{q,\la}^v: M'_{q,\la}\to M'_{q,\mu}\ten M'_{q,\ga}$
if $\Phi_{q,\la}^v(f^{\la+1}\cd v)=0$. This last condition is
satisfied if $\la+\mu\ge\ga$.

We conclude that, for $\la,\mu,\ga\in\Zplus$ and $\la-\mu-\ga$ even,
we can bring the intertwining map
$\Phi_{q,\la}^v: M_{q,\la}\to M_{q,\mu}\ten M_{q,\ga}$
down to the level of the corresponding finite dimensional modules
iff $|\ga-\mu|\le\la\le\ga+\mu$.
\end{remark}
We will later need the following observation, which immediately follows
from the existence and uniqueness of the intertwining map.
\begin{lemma} \label{th:IM27}
Let $V,\la,\mu,v$ and further assumptions be as in Lemma \ref{th:IM2}.
Let $W$ be another $\Uq$-module, and let $A: V\to W$ be
an $\Uq$-intertwining map. Then
\begin{equation}
(\id\ten A)\circ\Phi_{q,\la}^v=\Phi_{q,\la}^{Av}.
\label{eq:IM28}
\end{equation}
\end{lemma}
For a given $v \in V[\la -\mu]$, we want to determine 
the coefficients of the map $\Phi_{q,\la}^{v}$
on any element of the Verma module $M_{q,\la}$. 
Hence, we apply $f^{n}$ on both sides of 
\eqref{eq:IM4} and by use of the intertwining 
property of $\Phi_{q}$ and by
\eqref{{eq:PR6}}, we find the following expression:
\[
	\Phi_{q,\la}^{v}(f^{n}\cd x_{\la}) 
	=
	\sum_{j=0}^{n} \sum_{k=0}^{\iy} 
	{\frac{[n]_{q}!} {[j]_{q}!\, [n-j]_{q}!}}
	{\frac{q^{{\frac{1}{4}}(-\mu n + \la j)} 
	q^{-\frac{k}{4}( \la - 2n +2)}}{[k]_{q}!\,
	([-\mu]_{q})_{k}}}
	f^{k+j}\cd x_{\mu} \ten f^{n-j} e^{k}\cd v.
\]
Substitution of $\mu = \la - \wt(v)$ yields:
\begin{eqnarray*}
	\Phi_{q,\la}^{v}(f^{n}\cd x_{\la}) 
	= 
	\sum_{k=0}^{\iy} {\frac{ q^{-\frac{k}{4}
	(\la -2n +2)}}
	{[k]_{q}!\, ([-\la +  \wt (v)]_{q})_{k}}}
	\sum_{j=0}^{n} \frac{[n]_{q}!}{[j]_{q}!\, [n-j]_{q}!}
	q^{{\frac{1}{4}}(\wt (v) n -\la n + \la j)}&&\\
	\times f^{k+j} \cd x_{\la - \wt(v)} \ten f^{n-j} e^{k}\cd v.&&
\end{eqnarray*}
Finally, with new summation variables $m,j$, where $m= k+j$, we obtain:
\begin{equation}
	\Phi_{q,\la}^{v}(f^{n}\cd x_{\la}) 
	= 
	\sum_{m=0}^{\iy} f^{m}\cd x_{\la - \wt(v)} 
	\ten \FSF_{q,m,n}(\la)\cd v,
	\label{eq:IM6}
\end{equation}
where
\begin{equation}
\FSF_{q,m.n} (\la) 
= 
q^{\frac{m}{4}(2n - \la -2)-\frac{\la n}{4}}
\sum_{j=0}^{m \wedge n} q^{-\frac{j}{2}(n-\la-1)} 
{\frac{[n]_{q}!}{[j]_{q}! [m-j]_{q}!\,[n-j]_{q}!}}
f^{n-j}e^{m-j} {\frac{q^{{\frac{n}{4}}h}}{
([-\la +h]_{q})_{m-j}}}.
\label{eq:IM1}
\end{equation}
Then $\wt(\FSF_{q,m,n}(\la)\cd v)=\wt(v)+2m-2n$.
Note that the coefficients of the intertwining map are rational 
functions of $q^{\la}$.

Particular cases, which will be used to find the expression of
the fusion matrix \eqref{eq:FM1}, are the cases
$n=0$ and $m=0$:
\begin{equation}
	\FSF_{q,m,0}(\la) 
	= 
	{\frac{q^{-{\frac{m}{4}}(\la + 2)}}{[m]_{q}!}} e^{m} 
	{\frac {1}{([-\la + h]_{q})_{m}}}
	\label{eq:IM10}
\end{equation}
and
\begin{equation}
	\FSF_{q,0,n}(\la) 
	= 
	q^{-{\frac{\la n}{4}}} f^{n} q^{{\frac{n}{4}}h}.
	\label{eq:IM11}
\end{equation}

In the following, we want to find the expression of the intertwining 
operator in the particular case where the 
$\Uq$-module $V$ is a Verma module. 
Let $\la,~\mu,~\ga~\in \mathbb{C}$ with 
$\mu \notin \mathbb{Z}_{\geq 0}$ and 
$\mu + \ga -\la ~\in 2\mathbb{Z}_{\geq 0}$.
The intertwining operator  
$\Phi_{q,\la}^{x^{\ga}_{\la - \mu}}:
M_{q,\la} \longrightarrow M_{q,\mu} \ten M_{q,\ga}$ 
can then be written as
\begin{equation}
\Phi_{q,\la}^{x^{\ga}_{\la - \mu}} (x^{\la}_{\la - 2n}) 
= 
\sum_{m=0}^{n+{\frac{1}{2}}(\mu + \ga - \la)}
\FSC_{q,\mu - 2m, \la - \mu + 2m - 2n, \la - 2n}^{\mu, \ga, \la}\;
x_{\mu - 2m}^{\mu} 
\ten x_{\la - \mu + 2m - 2n}^\ga\,,
\label{eq:IM12}
\end{equation}
where the 
{\em generalized Clebsch-Gordan coefficients}
$\FSC_{q,\ldots}^{\;\cdots}$ satisfy
\[
\FSF_{q,m,n}(\la)\cd x_{\la -\mu}^\ga=
\FSC_{q,\mu - 2m, \la - \mu + 2m - 2n, \la - 2n}^{\mu, \ga, \la}\,
x_{\la - \mu + 2m - 2n}^\ga.
\]
Note that, by \eqref{eq:IM4}, we have
\begin{equation}
\FSC_{q,\mu, \la - \mu , \la}^{\mu, \ga, \la}=1.
\label{eq:IM24}
\end{equation}
\begin{theorem} \label{th:IM21}
The coefficients $\FSC_{q,\ldots}^{\;\cdots}$, defined by
\eqref{eq:IM12},
can be expressed 
in terms of $q$-Hahn polynomials \eqref{eq:PR40} or in terms of
$q$-hypergeometric functions as follows:
\begin{eqnarray}
&&		\FSC_{q,\mu-2m,\ga-2N+2m,\mu+\ga-2N}^
		{\mu,\ga,\mu+\ga-2l} 
		=
		q^{{\frac{m}{4}}(\mu + \ga)-{\frac{\mu}{4}}(N-l)-\half m(N-m)}
		\qbinom Nmq
		Q_l(q^{-m};q^{-\mu-1},q^{-\ga-1},N;q)\qquad\quad
		\label{eq:IM13}\\
&&\qquad\qquad\qquad\qquad
=q^{\frac{m}{4}(\mu + \ga)-\frac{\mu}{4}(N-l)-\half m(N-m)}
\qbinom Nmq	
\qhyp32{q^{-l},q^{-\ga-\mu+l-1},q^{-m}}{q^{-N},q^{-\mu}}{q,q}.
\label{eq:IM14}
	\end{eqnarray}
\end{theorem}
\Proof
By comparison of \eqref{eq:IM6} 
for $v= x^{\ga}_{ \ga - 2l}$
with \eqref{eq:IM12} and by the substitutions
$\la= \mu + \ga -2l$, $n= N-l$ 
($ N\in\Zplus$ and $0\le l\le N$) we get:
\[
	\FSF_{q,m, N-l} (\mu + \ga - 2l)\cd
	x^{\ga}_{\ga - 2l} 
	=
	 \FSC_{q,\mu - 2m, \ga - 2N + 2m, \mu+ \ga - 2N}^
	{\mu, \ga,\mu +\ga-2l}\, x^{\ga}_{\ga + 2m - 2N}, 
\]
where, in view of \eqref{eq:IM1}:
\begin{eqnarray*}
	&&\FSF_{q,m,N-l} (\mu + \ga - 2l) \cd
	x^{\ga}_{\ga - 2l}\\
	&&= 
	\sum_{j=0 \vee (m-l)}^{m \wedge (N-l)} 
	{\frac{q^{{\frac{m}{4}}(2N -\mu - \ga - 2)
	-{\frac{N-l}{4}}(\mu + \ga - 2l)-{\frac{j}{2}}
	(N+l-\mu - \ga -1)}
	[N-l]_{q}!}{[j]_{q}!\, [N-l-j]_{q}!\, [m-j]_{q}!}}\\
	&&\qquad\times
	 f^{N-l-j} e^{m-j} {\frac{q^{{\frac{N-l}{4}}h}}{
	([-\ga - \mu +2l +h]_{q})_{m-j}}} \cd
	x^{\ga}_{\ga-2l}\\ 
	&&=
	\sum_{j=0 \vee (m-l)}^{m \wedge(N-l)}
	{\frac{q^{{\frac{m}{4}}(2N-\mu - \ga - 2)
	-{\frac{\mu}{4}} (N-l)-{\frac{j}{2}}(N+l-\mu -\ga -1)}
	[N-l]_{q}!}{[j]_{q}!\,[N-l-j]_{q}!\,[m-j]_{q}!\,
	([-\mu]_{q})_{m-j}}}\,
	f^{N-l-j} e^{m-j}\cd x^{\ga}_{\ga - 2l}\\
	&&=
	\sum _{j=0 \vee (m-l)}^{m \wedge (N-l)} (-1)^{m-j}
	q^{{\frac{m}{4}}(2N-\mu - \ga - 2)-{\frac{\mu}{4}}(N-l)-
	{\frac{j}{2}}(N+l-\mu-\ga-1)}\\
	&&\qquad\times
	{\frac{\qfac{N-l}\, \qpoch{-l}{m-j} \qpoch{\ga-l+1}{m-j}}
	{\qfac j\, \qfac{N-l-j}\, \qfac{m-j}\, \qpoch{-\mu}{m-j}}} 
	\cd x^{\ga}_{\ga-2N +2m}.
\end{eqnarray*}
The last equality is obtained by use of \eqref{eq:PR13} and 
\eqref{eq:PR14}.
 Hence,
\begin{eqnarray}
	\FSC_{q,\mu-2m,\ga-2N+2m,\mu+\ga-2N}^{\mu,\ga,\mu+\ga-2l} 
	&=&
	\sum_{j=0\vee(m-l)}^{m \wedge (N-l)} (-1)^{m-j} 
	q^{{\frac{m}{4}}(2N-\mu-\ga-2) - {\frac{\mu}{4}}(N-l) 
	-{\frac{j}{2}(N+l-\mu -\ga-1)}}\nonumber\\
	&&\qquad\times
	{\frac{[N-l]_{q}!\, ([-l]_{q})_{m-j} ([\ga - l +1]_{q})_{m-j}}
	{[j]_{q}!\, [N-l-j]_{q}!\, [m-j]_{q}!\, ([-\mu]_{q})_{m-j}}} 
	\label{eq:IM15}.
\end{eqnarray}

In order to express the coefficients $\FSC_{q}$ in terms of $q$-Hahn 
polynomials, we separately study the cases $ m \geq l$ and $ m\le l$. 
In the first case ($m\geq l$), 
we change the summation variable to $i:=j-m+l$.
Then \eqref{eq:IM15}
can be rewritten as
\begin{eqnarray*}
	\FSC_{q,\mu - 2m, \ga - 2N + 2m, \mu +\ga -2N}^{\mu, \ga,\mu + 
	\ga -2l} 
	&=& 
	\sum_{i=0}^{(N-m) \wedge l}  q^{\frac{m}{4}( \mu + \ga - 2l) 
	-\frac{\mu}{4}(N-l)+\frac{l-i}{2}(N+l-\ga-\mu-1)}
	{\frac{(-1)^{l}[N-l]_{q}!\, ([-\ga]_{q})_{l} }{[m-l]_{q}!\, [N-m]_{q}!\,
	([-\mu]_{q})_{l}}}\\
	&&\qquad\qquad\times
	{\frac{ ([-l]_{q})_{i} ([-N+m]_{q})_{i} ([\mu - l +1]_{q})_{i}}
	{[i]_{q}!\, ([-\ga]_{q})_{i} ([m-l+1]_{q})_{i}}}.
\end{eqnarray*}
We have
\[
	{\frac{ ([-l]_{q})_{i} ([-N+m]_{q})_{i} ([\mu - l +1]_{q})_{i}}
	{[i]_{q}!\, ([-\ga]_{q})_{i} ([m-l+1]_{q})_{i}}}
	= 
	q^{-{\frac{i}{2}}(-N-l +\mu +\ga-1)} 
	{\frac{(q^{-l};q)_{i}(q^{-N+m};q)_{i} (q^{\mu -l +1};q)_{i}}
	{(q;q)_{i} (q^{-\ga};q)_{i} (q^{m-l+1};q)_{i}}}
\]
and
\begin{eqnarray*}
	{\frac{(-1)^{l} [N-l]_{q}!\, ([-\ga]_{q})_{l}}
	{[m-l]_{q}!\, [N-m]_{q}!\, ([-\mu]_{q})_{l}}}
	&=&
	\frac{[N]_{q}!}{[m]_{q}!\,[N-m]_{q}!} 
	\frac{([-\ga]_{q})_{l} ([m-l+1]_{q})_{l}}
	{([\mu-l+1]_{q})_{l}([N-l+1]_{q})_{l}}\\
	&=&
	q^{-\frac{l}{2}(-\ga - \mu + m - N +l -1)-\half m(N-m)}
\qbinom Nmq
	\frac{(q^{-\ga};q)_{l} (q^{m-l+1};q)_{l}}
	{(q^{\mu - l + 1};q)_{l} (q^{N-l+1};q)_{l}}.
\end{eqnarray*}
Hence,
\begin{eqnarray}
&&	\FSC_{q,\mu - 2m, \ga - 2N + 2m, \mu +\ga -2N}^
	{\mu, \ga,\mu + \ga -2l}	
	= 
	q^{\frac{m}{4}(\mu + \ga)-\frac{\mu}{4}(N-l)-\half m(N-m)}
	\qbinom Nmq	\nonumber
\\
&&\qquad\qquad\qquad\qquad\times
q^{l(N-m)}\frac{(q^{-\ga};q)_{l} (q^{m-l+1};q)_{l}}
{(q^{\mu - l + 1};q)_{l} (q^{N-l+1};q)_{l}}
\qhyp32{q^{-l}, q^{-N+m}, q^{\mu-l+1}}{q^{-\ga},q^{m-l+1}}{q,q}
\label{eq:IM25}
\\
&&\quad=q^{\frac{m}{4}(\mu + \ga)-\frac{\mu}{4}(N-l)-\half m(N-m)}
\qbinom Nmq	
\frac{(q^{-\ga};q)_{l}}{(q^{\mu - l + 1};q)_{l}}
\qhyp32{q^{-l},q^{-\ga-\mu+l-1},q^{-N+m}}{q^{-\ga},q^{-N}}{q,q^{\mu-m}}
\nonumber
\\
&&\quad=q^{\frac{m}{4}(\mu + \ga)-\frac{\mu}{4}(N-l)-\half m(N-m)}
\qbinom Nmq	
\qhyp32{q^{-l},q^{-\ga-\mu+l-1},q^{-m}}{q^{-N},q^{-\mu}}{q,q},
\nonumber
\end{eqnarray}
which is \eqref{eq:IM14} for $m\ge l$.
Above we used the transformation formula \eqref{eq:PR30} twice.
In these inequalities the occurrence of $q^{-N}$ as a denominator parameter
is not harmful for the application of \eqref{eq:PR30},
since all parts have $q^{-l}$ as a numerator parameter, while $l\le N$.

In the second case ($m \le l$), formula \eqref{eq:IM15} can be 
rewritten as
\begin{eqnarray*}
	\FSC_{q,\mu-2m,\ga-2N+2m,\mu+\ga-2N}^
	{\mu,\ga,\mu+\ga-2l}	
	&=&
	 q^{{\frac{m}{4}}(2N-\mu-\ga-2) -\frac{\mu}{4}(N-l) - 
	\frac{j}{2}(N+l-\mu-\ga-1)}
	{\frac{(-1)^{m}([-l]_{q})_{m}([\ga-l+1]_{q})_{m}}{[m]_{q}!\,
	([-\mu]_{q})_{m}}}\\
	&& \times
	 \sum_{j=0}^{m\wedge(N-l)} {\frac{([-N+l]_{q})_{j} ([-m]_{q})_{j} 
	([\mu -m +1]_{q})_{j}}
	{[j]_{q}!\, ([-\ga + l -m]_{q})_{j} ([l-m+1]_{q})_{j}}}.
\end{eqnarray*}
We have
\[
	{\frac{([-N+l]_{q})_{j} ([-m]_{q})_{j} 
	([\mu -m +1]_{q})_{j}}
	{[j]_{q}!\, ([-\ga + l -m]_{q})_{j} ([l-m+1]_{q})_{j}}}
	=
	q^{-\frac{j}{2}(\mu +\ga - N-l -1)}
	\frac{(q^{-N+l};q)_{j} (q^{-m};q)_{j} (q^{\mu - m + 1};q)_{j}}
	{(q;q)_{j} (q^{-\ga + l - m};q)_{j} (q^{l-m + 1};q)_{j}}
\]
and
\[
	\frac{(-1)^{m} ([-l]_{q})_{m} ([\ga - l +1]_{q})_{m}}
	{[m]_{q}!\, ([-\mu]_{q})_{m}}
	=
	q^{-\frac{m}{2}(-\mu - \ga - N + 2l - 1)-\half m(N-m)}
	\qbinom Nmq
	\frac{(q^{l-m+1};q)_{m} (q^{-\ga+l-m};q)_{m}}
	{(q^{N-m+1};q)_{m} (q^{\mu - m + 1};q)_{m}}.
\]
Hence,
\begin{eqnarray}
&&\FSC_{q,\mu-2m,\ga-2N+2m,\mu+\ga-2N}^{\mu,\ga,\mu+\ga-2l}
=
q^{\frac{m}{4}(\mu + \ga)-\frac{\mu}{4}(N-l)-\half m(N-m)}
\qbinom Nmq \nonumber\\
&&\qquad\qquad\times
q^{m(N-l)} {\frac{(q^{l-m+1};q)_{m}(q^{-\ga+l-m};q)_{m}}
{(q^{N-m+1};q)_{m}(q^{\mu-m+1};q)_{m}}}
\qhyp32{q^{-N+l},q^{\mu-m+1}, q^{-m}}{q^{l-m-\ga}, q^{l-m+1}}{q,q}
\label{eq:IM26}\\
&&\quad=q^{\frac{m}{4}(\mu + \ga)-\frac{\mu}{4}(N-l)-\half m(N-m)}\qbinom Nmq
\frac{(q^{-\ga+l-m};q)_m}{(q^{\mu-m+1};q)_m}
\qhyp32{q^{-m},q^{-\mu-\ga+l-1},q^{-N+l}}{q^{-N},q^{-\ga+l-m}}
{q,q^{\mu-l+1}} \nonumber\\
&&\quad=q^{\frac{m}{4}(\mu + \ga)-\frac{\mu}{4}(N-l)-\half m(N-m)}\qbinom Nmq
\qhyp32{q^{-m},q^{-\mu-\ga+l-1},q^{-l}}{q^{-N},q^{-\mu}}
{q,q}, \nonumber
\end{eqnarray}
which is \eqref{eq:IM14} for $m\le l$.
Above we used again the transformation formula \eqref{eq:PR30} twice.
$\square$\\

Let $\mu,\ga\in\CC\backslash\Zplus$ and
$l\in\Zplus$\,. Write $\la:=\mu+\ga-2l$.
By combination of the $\Uq$-intertwining properties of the mappings
\[
\Phi_{q,\la}^{f^l\cd x_\ga}: M_{q,\la}\to M_{q,\mu}\ten M_{q,\ga},
\quad
\Phi_{q,\la}^{f^l\cd x_\mu}: M_{q,\la}\to M_{q,\ga}\ten M_{q,\mu},
\quad
P\circ\FSR: M_{q,\mu}\ten M_{q,\ga}\to M_{q,\ga}\ten M_{q,\mu}
\]
we can expect a nice relationship between the operators
$P\circ\FSR\circ\Phi_{q,\la}^{f^l\cd x_\ga}$ and
$\Phi_{q,\la}^{f^l\cd x_\mu}$.
In fact, we have:
\begin{theorem}
\label{th:IM17}
Let $\mu,\ga\in\CC\backslash\Zplus$ and $l\in\Zplus$. Then
\begin{equation}
P\circ\FSR\circ\Phi_{q,\mu+\ga-2l}^{f^l\cd x_\ga}
=
(-1)^lq^{-\quart l(\mu+\ga-2l+2)+\quart\ga(\mu-2l)}\,
\frac{\qpoch{-\ga}l}{\qpoch{-\mu}l}\,
\Phi_{q,\mu+\ga-2l}^{f^l\cd x_\mu}\,.
\label{eq:IM16}
\end{equation}
\end{theorem}
\Proof
Combination of \eqref{eq:IM6} with \eqref{eq:PR7} yields
\begin{eqnarray}
&&P\left(\FSR\cd\Phi_{q,\mu+\ga-2l}^{f^l\cd x_\ga}(x_{\mu+\ga-2l})\right)
\nonumber\\
&&=
\sum_{k=0}^l \sum_{j=0}^k
\frac{(1-q^{-1})^j q^{-\quart j(j-1)} q^{-\quart k(\mu+\ga-2l+2)}}
{\qfac j\, \qfac k\, \qpoch{-\mu}k}\,
q^{\quart(h\ten h)}
\left(q^{-\quart jh} f^j e^k f^l\cd x_\ga\ten q^{\quart jh} e^j
f^k\cd x_\mu\right)\qquad\quad
\label{eq:IM17}\\
&&=
\frac{q^{-\quart l(\mu+\ga-2l+2)} q^{\quart\ga(\mu-2l)}}
{\qfac l\,\qpoch{-\mu}l}\,e^lf^l\cd x_\ga\ten f^l\cd x_\mu
+M_{q,\ga}[\wt<\ga]\ten M_{q,\mu}\nonumber\\
&&=(-1)^l q^{-\quart l(\mu+\ga-2l+2)+\quart\ga(\mu-2l)}\,
\frac{\qpoch{-\ga}l}{\qpoch{-\mu}l}\,x_\ga\ten f^l\cd x_\mu
+M_{q,\ga}[\wt<\ga]\ten M_{q,\mu}.\nonumber
\end{eqnarray}
Finally combine the intertwining property of $P\circ\FSR$ with Definition
\ref{th:IM5}. $\square$\\

Of course, an independent verification of \eqref{eq:IM16} should be
possible without use of the intertwinining property of $P\circ\FSR$.
For that purpose, simplify the double sum \eqref{eq:IM17},
replace the summation variable $j$ by a summation variable
$i:=l-k+j$, write the resulting double sum as $\sum_{i=0}^l\sum_{j=0}^i$\,,
and reduce the inner sum to a ${}_2\phi_0$ sum which can be evaluated by use
of \eqref{eq:PR29}:
${}_2\phi_0(q^{-i},q^{\ga-i+1};-;q,q^{2i-\ga-1})=q^{i(i-\ga-1)}$.
\section{Orthogonality of q-Hahn polynomials and the\\
Shapovalov form} \label{sec:OS}
In this section we will derive the known orthogonality relations
\begin{eqnarray}
	&&\sum_{m=0}^{N}\frac{(q^{-\mu},q^{-N};q)_{m}}
	{(q,q^{\ga-N+1};q)_{m}}\, q^{m(\mu + \ga + 1)}\,
	Q_l(q^{-m};q^{-\mu-1},q^{-\ga-1},N;q)\,
	Q_{l'}(q^{-m};q^{-\mu-1},q^{-\ga-1},N;q)\qquad
\nonumber\\
	&&\quad=\de_{l,l'}
	\frac{q^{\mu N}(q^{-\mu-\ga};q)_{N}}{(q^{-\ga};q)_{N}}\,
	\frac {(q,q^{-\mu - \ga +N},q^{-\ga};q)_{l}}
	{(q^{-\mu}, q^{-\mu - \ga - 1}, q^{-N} ; q)_{l}}\,
\frac{1-q^{-\mu-\ga-1}}{1- q^{-\mu - \ga + 2l -1}}\,
(-1)^l q^{\half l(l-1)}q^{-(N+\mu)l}
	\label{eq:OS10}
\end{eqnarray}
for $q$-Hahn polynomials (see \cite[(3.6.2)]{KoeS})
as a consequence of the quantum group interpretation
\eqref{eq:IM13} of $q$-Hahn prolynomials.
As a tool we will use the quantum analogue of the so-called
Shapovalov form, see for instance \cite[\S5]{DoT} and references given there.

let $V$ be a $\Uq$-module on which a symmetric
bilinear form $\lan\;\ran$ is given. We will call this form
a {\em (quantum) Shapovalov form} on $V$ if
\begin{equation}
\lan e\cd v,w\ran=\lan v,f\cd w\ran
\quad\mbox{and}\quad
\lan h\cd v,w\ran=\lan v,h\cd w\ran
\label{eq:OS1}
\end{equation}
for all $v,w\in V$.
We observe the following properties of the Shapovalov form.

On the Verma module $M_{q,\la}$ we find by
\eqref{eq:PR13}, \eqref{eq:PR14} that a Shapovalov form exists uniquely,
up to a constant factor:
\begin{equation}
	\lan f^{n}\cd x_{\la}, f^{m}\cd x_{\la}\ran  = \de_{m,n} 
	(-1)^{n} [n]_{q}!\, 
	([-\la]_{q})_{n}\,\lan x_\la,x_\la\ran.
	 \label{eq:OS11}
\end{equation}
By convention we normalize the Shapovalov form on $M_{q,\la}$ by putting
$\lan x_\la,x_\la\ran:=1$.

Next, let $V$ and $W$ be $\Uq$-modules equipped with a Shapovalov form.
Then, in view of \eqref{eq:PR4}, we can define a Shapovalov form on 
$V\ten W$ by
\begin{equation}
\lan v_1\ten w_1,v_2\ten w_2\ran=\lan v_1,v_2\ran\,\lan w_1,w_2\ran.
\label{eq:OS12}
\end{equation}

As a final property, let $V$ be a $\Uq$-module equipped with a
Shapovalov form and with
highest weight submodules $W,W'$ with
highest weights $\la$ resp.\ $\mu$ such that $\la\ne\mu$ and $\la\ne -\mu-2$.
Then $\lan w,w'\ran=0$ for $w\in W$, $w'\in W'$, because,
by \eqref{eq:PR15}, we have
\[
[\thalf\la]_q[\thalf\la+1]_q\lan w,w'\ran=\lan\Om\cd w,w'\ran=
\lan w,\Om\cd w'\ran=[\thalf\mu]_q[\thalf\mu+1]_q\lan w,w'\ran.
\]

For a given Shapovalov form $\lan\;\ran$ on a $\Uq$-module $V$
we will also use the notation
\begin{equation}
\|v\|^2:=\lan v,v\ran\qquad(v\in V).
\label{eq:OS16}
\end{equation}
If $\|v\|^2>0$ then we will also work with $\|v\|$, being the positive
square root of $\|v\|^2$.
\begin{lemma}
Let $\mu,\ga\in\CC$, $n\in\Zplus$. Assume that $\mu\notin\Zplus$.
Then
\begin{equation}
\|\Phi_{q,\mu + \ga - 2l}^{f^{l}\cd x_{\ga}}
	(x_{\mu + \ga - 2l})\|^2
=(-1)^l q^{\half l(l-\ga-1)} 
	 \frac{\qfac l\,\qpoch{-\ga}l\qpoch{-\mu-\ga+l-1}l}
	{\qpoch{-\mu}l}\,.
\label{eq:OS13}
\end{equation}
\end{lemma}
\Proof
By \eqref{eq:IM4}, \eqref{eq:OS12},
\eqref{eq:OS11} and \eqref{eq:PR13} we obtain:
\begin{eqnarray*}
	 &&\lan \Phi_{q,\mu + \ga - 2l}^{f^{l}\cd x_{\ga}}
	(x_{\mu + \ga - 2l}),\Phi_{q,\mu +
\ga - 2l}^{f^{l}\cd x_{\ga}} (x_{\mu + \ga - 2l}) \ran \\
	&&\qquad\qquad =
	 \sum_{k=0}^{\iy}
 	q^{-\frac{k}{2}(\mu + \ga - 2l +2)} \frac{1} 
	{\left( [k]_{q}!\,([-\mu]_{q})_{k} \right)^{2}}
	\lan f^{k}\cd x_{\mu}, f^{k}\cd x_{\mu} \ran
	\lan e^{k} f^{l} x_{\ga},e^{k} f^{l} x_{\ga} \ran\\
&&\qquad\qquad=(-1)^{l}[l]_{q}!\, ([-\ga]_{q})_{l}
\sum_{k=0}^l\frac{([-l]_q)_k([\ga-l+1]_q)_k}{([-\mu]_q)_k[k]_q!}
q^{-\half k(\mu+\ga-2l+2)}\\
&&\qquad\qquad=(-1)^{l}[l]_{q}!\, ([-\ga]_{q})_{l}\, 
	\qhyp21{q^{-l},q^{\ga - l +1}}{q^{-\mu}}{q,q^{-(\mu + \ga - 2l +1)}}.
\end{eqnarray*}
Now \eqref{eq:OS13} follows by \eqref{eq:PR28}
(the reversed $q$-Chu-Vandermonde sum).
$\square$\\[\bigskipamount]\noindent
{\bf Proof of (\ref{eq:OS10})}\quad
By continuity it is sufficient to prove \eqref{eq:OS10}
if moreover $\mu+\ga\notin\ZZ$.
If $l \ne l'$ then
$\Phi_{q,\mu + \ga - 2l}^{f^{l}\cd x_{\ga}}(M_{q,\mu+\ga-2l})$ and
$\Phi_{q,\mu + \ga - 2l'}^{f^{l'}\cd x_{\ga}}(M_{q,\mu+\ga-2l'})$
are highest weight submodules of $M_{q,\mu}\ten M_{q,\ga}$ with
distinct highest weights $\mu+\ga-2l$ and $\mu+\ga-2l'$, respectively.
Hence, the Shapovalov form on $M_{q,\mu}\ten M_{q,\ga}$ with
its two arguments restricted to
$\Phi_{q,\mu + \ga - 2l}^{f^{l}\cd x_{\ga}}(M_{q,\mu+\ga-2l})$ and
$\Phi_{q,\mu + \ga - 2l'}^{f^{l'}\cd x_{\ga}}(M_{q,\mu+\ga-2l'})$,
respectively,
yields 0. Thus
\begin{equation}
\lan\Phi_{q,\mu + \ga - 2l}^{f^{l}\cd x_{\ga}}
(f^{N-l}\cd x_{\mu + \ga - 2l}),
\Phi_{q,\mu + \ga - 2l'}^{f^{l'}\cd x_{\ga}}
(f^{N-l'}\cd x_{\mu + \ga - 2l'})\ran=0\qquad(l\ne l').
\label{eq:OS14}
\end{equation}
For $l=l'$ we obtain by the intertwining property of $\Phi_{q}$ and by
\eqref{eq:PR13}, \eqref{eq:OS1} and \eqref{eq:OS13}
that
\begin{eqnarray}
&&\nonumber\lan\Phi_{q,\mu + \ga - 2l}^{f^{l}\cd x_{\ga}}
(f^{N-l}\cd x_{\mu + \ga - 2l}),
\Phi_{q,\mu + \ga - 2l'}^{f^{l'}\cd x_{\ga}}
(f^{N-l'}\cd x_{\mu + \ga - 2l'})\ran\\
&&\quad\nonumber=
	[N-l]_{q}!\, ([-\mu - \ga + 2l]_{q})_{N-l} (-1)^{N-l} 
	\lan \Phi_{q,\mu + \ga - 2l}^{f^{l}\cd x_{\ga}}
	(x_{\mu + \ga - 2l}),
	\Phi_{q,\mu + \ga - 2l}^{f^{l}\cd x_{\ga}}
	(x_{\mu + \ga - 2l}) \ran\\
&&\quad=\frac{[N]_{q}!\, ([-\mu -\ga]_{q})_{N} 
	[-\mu -\ga - 1]_{q}}{[-\mu -\ga + 2l -1]_{q}}
	\frac{(-1)^{N+l} q^{\frac{l}{2}(l-\ga-1)} [l]_{q}!\, 
	([-\ga]_{q})_{l} ([-\mu - \ga + N ]_{q})_{l}}
	{ ([-\mu]_{q})_{l} ([-N]_{q})_{l} 
	([-\mu - \ga - 1 ]_{q})_{l}}.\qquad
\label{eq:OS15}
\end{eqnarray}
On the other hand, by use of
\eqref{eq:IM12}, \eqref{eq:IM13} we can write:
\begin{eqnarray*}
	&&\lan \Phi_{q,\mu + \ga - 2l}^{f^{l}\cd x_{\ga}}
	(f^{N-l}\cd x_{\mu + \ga - 2l}),
	\Phi_{q,\mu + \ga - 2l'}^{f^{l'}\cd x_{\ga}}
	(f^{N-l'}\cd x_{\mu + \ga - 2l}) \ran \\
	&&=
	\sum_{m=0}^{N} \lan \FSC_{q, \mu - 2m, \ga - 2N +2m,
	 \mu + \ga -2N}^{\mu, \ga,\mu + \ga - 2l}
	x_{\mu - 2m}^{\mu} \ten x_{\ga - 2 N + 2m }^{ \ga}, 
 	\FSC_{q, \mu - 2m, \ga - 2N +2m, \mu + \ga -2N}^
	{\mu, \ga, \mu + \ga - 2l'} x_{\mu - 2m}^{\mu} 
	\ten x_{\ga - 2 N + 2m }^{ \ga} \ran\\
	&&=
	\sum_{m=0}^{N} 
	\left( \frac{[N]_{q}!}{[m]_{q}!\, [N-m]_{q}!}\right)^{2}
	 q^{\frac{m}{2}(\mu + \ga)-\frac{\mu}{4}(2N-l-l')}
	(Q_{l}Q_{l'})(q^{-m}, q^{-\mu - 1}, q^{-\ga - 1}, N;q) \\
	&&\qquad\qquad\times\lan f^{m}\cd x_{\mu}, f^{m}\cd x_{\mu}\ran
	\lan f^{N-m}\cd x_{\ga},f^{N-m}\cd x_{\ga}\ran \\
	&&=
	(-1)^{N} \left( [N]_{q}!\, \right)^{2}  \sum_{m=0}^{N}
	q^{\frac{m}{2}(\mu + \ga)-\frac{\mu}{4}(N-l-l')}
	\frac{([-\mu]_{q})_{m}([-\ga]_{q})_{N-m}}{[m]_{q}!\, [N-m]_{q}!}
	(Q_{l}Q_{l'})(q^{-m}, q^{-\mu - 1}, q^{-\ga - 1}, N;q).	
\end{eqnarray*}
Combination of this last result with \eqref{eq:OS15} and
\eqref{eq:OS14} yields \eqref{eq:OS10}. $\square$

\begin{remark} \label{th:OS2}
Use the notation of Remark \ref{th:PR20}.
Let $j_1,j_2\in\thalf\Zplus$. Then the tensor product $V^{j_1}\ten V^{j_2}$
of two finite dimensional irreducible $\Uq$-modules decomposes as the direct
sum of all $V^j$ such that $j\in\{|j_1-j_2|, |j_1-j_2|+1,\ldots,j_1+j_2\}$.
The Shapovalov forms on $V^{j_1}$ and $V^{j_2}$ induce a Shapovalov
form on $V^{j_1}\ten V^{j_2}$ and hence on each $V^j$ occurring in this
tensor product.
Let the vectors $e_m^j(j_1,j_2)$ ($m=-j,-j+1,\ldots,j$)
form the standard basis of the irreducible submodule $V^j$
of $V^{j_1}\ten V^{j_2}$, where the basis vectors are orthonormal
with respect to the Shapovalov form.
This basis is unique up to a constant complex
factor of absolute value 1, independent of $j$. Normalize the basis such
that the inner product between $e_j^j(j_1,j_2)$ and
$e_{j_1}^{j_1}\ten e_{j-j_1}^{j_2}$ is positive.
There will be an expansion of the form
\begin{equation}
e_m^j(j_1,j_2)=\sum_{m_1+m_2=m}
\qthreej{j_1}{j_2}j{m_1}{m_2}m
e_{m_1}^{j_1}\ten e_{m_2}^{j_2}.
\label{eq:OS3}
\end{equation}
The coefficients in \eqref{eq:OS3} are called {\em $q$-$3j$ symbols} or
{\em $q$-Clebsch-Gordan coefficients}. See for instance
\cite{KiR} for further discussion.

Because of Remark \ref{th:IM9} formula \eqref{eq:IM12} remains valid
by analytic continuation for the case of finite dimensional irreducible
$\Uq$-modules with $\la,\mu,\ga$ as at the end of Remark \ref{th:IM9}.
In that case, In view of Definition \ref{th:IM5}, we can relate the
$q$-Clebsch-Gordan coefficients as defined by \eqref{eq:IM12} to
the $q$-$3j$ symbols in \eqref{eq:OS3}.
First we specialize \eqref{eq:IM12} and \eqref{eq:OS13}
to the finite dimensional case:
\begin{equation}
\Phi_{q,2j}^{x_{2j-2j_1}^{2j_2}}(x_{2m}^{2j})=
\sum_{m_1+m_2=m}\FSC_{q,2m_1,2m_2,2m}^{2j_1,2j_2,2j}\;
x_{2m_1}^{2j_1}\ten x_{2m_2}^{2j_2},
\label{eq:OS4}
\end{equation}
\begin{eqnarray}
&&\|\Phi_{q,2j}^{x_{2j-2j_1}^{2j_2}}(x_{2j})\|^2
=q^{\half(j_1+j_2-j)(j_1-j_2-j-1)}
\nonumber\\
&&\qquad\qquad\qquad\quad
\times
\frac{\qfac{j_1+j_2-j}\,\qpoch{j_2-j_1+j+1}{j_1+j_2-j}\,
\qpoch{2j+2}{j_1+j_2-j}}
{\qpoch{j_1-j_2+j+1}{j_1+j_2-j}}\,.
\label{eq:OS5}
\end{eqnarray}
Note that $\FSC_{q,2j-2j_2,2j_2,2j}^{2j_1,2j_2,2j}=1$ by \eqref{eq:IM24},
and that the \RHS\ of \eqref{eq:OS5} is $>0$.

The normalized intertwining operator
\begin{equation}
\wtilde\Phi_{q,j}^{e^{j_2}_{j-j_1}}:=
\|\Phi_{q,2j}^{e_{j-j_1}^{j_2}}(e_j^j)\|^{-1}\,
\Phi_{q,2j}^{e^{j_2}_{j-j_1}}=
\|\Phi_{q,2j}^{x_{2j-2j_1}^{2j_2}}(x_{2j})\|^{-1}\,
\Phi_{q,2j}^{x_{2j-2j_1}^{2j_2}}
\label{eq:OS6}
\end{equation}
will thus satisfy
\begin{equation}
	\wtilde{\Phi}_{q,j}^{e^{j_2}_{j-j_1}}(e_m^j)=e_m^j(j_1,j_2).
\label{eq:OS7}
\end{equation}
We conclude from \eqref{eq:OS4}, \eqref{eq:OS5},
\eqref{eq:OS6}, \eqref{eq:OS7}, \eqref{eq:OS3} and \eqref{eq:PR21} that
\begin{eqnarray}
	&&\qthreej{j_1}{j_2}j{m_1}{m_2}m
	=
\left(\frac{q^{-\half (j_1 +j_2 -j)(j_1 - j_2 -j -1)}
\qpoch{j_1-j_2+j+1}{j_1 + j_2 -j}}
{\qpoch{2j+2}{j_1+j_2-j}\,\qpoch{j_2-j_1+j+1}{j_1+	_2-j}\qfac{j_1 +j_2 -j}}
\right)^\half
\nonumber\\
&&\qquad\qquad\qquad\times \left(
\frac
{\qpoch{j_1+m_1+1}{j_1-m_1}\,\qfac{j_1-m_1}\,
\qpoch{j_2+m_2+1}{j_2-m_2}\,\qfac{j_2 - m_2}}
{\qpoch{j+m+1}{j-m}\,\qfac{j-m}}\right)^\half
\nonumber\\
&&\qquad\qquad\qquad\times\FSC_{q, 2m_1,2m_2,2m}^{2j_1,2j_2,2j}
\qquad\qquad\qquad\qquad\qquad\qquad\qquad\quad
(m_1+m_2=m).
\label{eq:OS8}
\end{eqnarray}

Specialization of \eqref{eq:IM15} to the finite dimensional case,
in combination with \eqref{eq:OS8},  yields the $q$-analogue
\cite[(3.4)]{KiR} of the Racah formula for Clebsch-Gordan coefficients.
(However, note that in \cite[(3.4)]{KiR} the factor $q^{\half m_1(m_1+1)}$
is not correct. It should be replaced by $q^{\half m(m_1+1)}$.)$\;$

As we saw in the proof of Theorem \ref{th:IM21}, the reduction of
\eqref{eq:IM15} to the final $q$-hypergeometric form \eqref{eq:IM14}
passed through
an intermediate $q$-hypergeometric form occurring in two versions
\eqref{eq:IM25} and \eqref{eq:IM26} depending on the sign of $m-l$.
Specialization of \eqref{eq:IM26} to the finite dimensional case,
in combination with \eqref{eq:OS8}, yields for $j-j_2-m_1\le0$ the
$q$-analogue
\cite[(3.60)]{BiLoh}) of the Racah formula for Clebsch-Gordan coefficients.
\end{remark}

In the following we prove the following result stated in
\cite[(4.8)]{KiR}. 
\begin{equation}
P\FSR e_m^j(j_1,j_2)=(-1)^{j-j_1-j_2}~ q^{\half(c_j - c_{j_1} -c_{j_2})}
e_m^j(j_2,j_1),
\label{eq:AP2}
\end{equation}
where $c_j:=j(j+1)$.
Graphically this formula reads as follows.
%
\vspace{.8cm}

\hspace{2.2cm}$P\cal{R}$

\vspace{-1.2cm}

\hspace{3cm}\includegraphics[height=2cm]{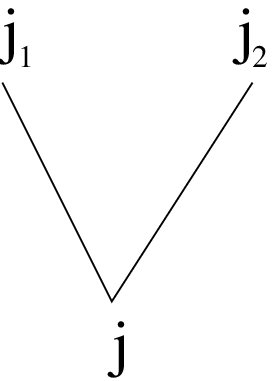}

\vspace{-1.3cm}

\hspace{5cm}$=~ (-1)^{j-j_1-j_2}~ q^{\half(c_j - c_{j_1} -c_{j_2})}$

\vspace{-1.4cm}

\hspace{10cm}\includegraphics[height=2cm]{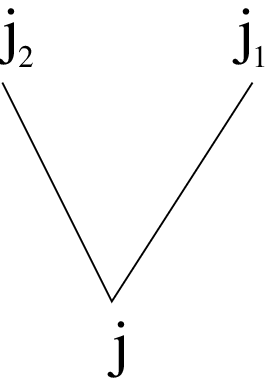}

\vspace{-1cm}

\begin{equation}
\label{eq:AP3}
\end{equation}
Here, for the diagram both on the left and on the right one has to
substitute the $m$th standard basis vector in the module $V_j$, which
is determined  within the tensor product by the diagram.
\\[\medskipamount]
{\bf Proof of (\ref{eq:AP2})}\quad
Formula \eqref{eq:AP2} 
can be equivalently written in terms of normalized intertwining
operators \eqref{eq:OS6} as follows:
\begin{equation}
P\circ\FSR\circ\wtilde{\Phi}_{q,j}^{e_{j-j_1}^{j_2}} =
(-1)^{j-j_1-j_2}~ q^{\half(c_j - c_{j_1} -c_{j_2})}\,
\wtilde{\Phi}_{q,j}^{e_{j-j_2}^{j_1}}.
\label{eq:AP4}
\end{equation}
We will obtain \eqref{eq:AP4} from \eqref{eq:IM16}, which remains
valid for finite dimensional relations in view of
Remark \ref{th:OS2}, and which can then be written as
\begin{eqnarray}
	&& P\circ\FSR\circ\Phi_{q,2j}^{ x_{2(j-j_1)}^{2j_2}}=
	(-1)^{j_1+j_2-j}q^{-\half (j_1+j_2-j)(j+1)+j_2(j-j_2)}
\nonumber\\
&&\qquad\qquad\qquad\qquad\qquad\times\
\frac{\qpoch{j_2-j_1+j+1}{j_1+j_2-j}}{\qpoch{j_1-j_2+j+1}{j_1+j_2-j}}\,
\Phi_{q,2j}^{x_{2(j-j_2)}^{2j_1}}\label{app3}.
\label{eq:AP5}
\end{eqnarray}
Formula \eqref{eq:AP4} now follows from \eqref{eq:AP5} by use of
\eqref{eq:OS6} and \eqref{eq:OS5}.
$\square$

\section{The fusion matrix} \label{sec:FM}
Let $v\in V$ and $w\in W$ be weight vectors in
$\Uq$-modules $V$ and $W$ and 
let $\la\in\CC$ such that
$\la - \wt(v),\;\la - \wt(v) - \wt(w)\notin\Zplus$.
Then the composition map
\[
	\Phi_{q,\la}^{w, v}: M_{q,\la } 
	\stackrel{\Phi_{q,\la}^{v}} {\longrightarrow}
	M_{q,\la - \wt(v)} \ten V 
	\stackrel{\Phi_{q,\la-\wt(v)}^{w} \ten \id}
	{\longrightarrow} M_{q,\la - \wt(v) - \wt(w)} \ten W \ten V
\]
is $\Uq$-intertwining. By the considerations of Section 3
this composition map must have the form $\Phi_{q,\la}^u$
for some $u\in(W\ten V)[\wt(v)+\wt(w)]$. It will turn out that
the map $w\ten v\mapsto u$ can be linearly extended to a map
$J_{W,V}(\la): W\ten V\to W\ten V$. We call this operator
the {\em fusion matrix} for $W$ and $V$, see
\cite[\S2.1]{ES}. Note its defining property
\begin{equation}
(\Phi_{q,\la-\wt(v)}^w\ten\id)\circ\Phi_{q,\la}^v=
\Phi_{q,\la}^{J_{W,V}(\la)(w\ten v)}.
\label{eq:FM1}
\end{equation}
\begin{theorem}
The fusion matrix can be written as:
\begin{equation}
	J_{W,V}(\la) (w \ten v)
	= 
	\sum_{l= 0}^{\iy } 
	\frac{q^{-{\frac{l}{2}}(\la+1)}}{[l]_{q}!}
f^{l} q^{lh/4}\cd w
	\ten
e^{l}\frac{q^{lh/4}}{([-\la + h]_{q})_{l}}\cd v.
	\label{eq:FM2}
\end{equation}
\label{th:FM3}
\end{theorem}
\Proof
We first compute the composition map of the intertwining operators
by use of \eqref{eq:IM6}:
\begin{eqnarray*}
	\Phi_{q,\la}^{w,v} (x_{\la})\hskip-0.5truecm
	&&= 
	(\Phi^{w}_{q,\la - \wt(v)} \ten \id) 
	\circ \Phi_{q,\la}^{v} (x_{\la})\\
	&&=
	(\Phi^{w}_{q,\la - \wt(v)} \ten\id) \sum_{m=0}^{\iy}
	f^{m} \cd x_{\la - \wt (v)} \ten \FSF_{q,m,0}(\la)\cd v\\
	&&= 
	\sum_{m=0}^{\iy}\sum_{m'=0}^{\iy} 
	f^{m'}\cd x_{\la - \wt(v) - \wt( w)} \ten 
	\FSF_{q,m',m}(\la - \wt(v))\cd w \ten \FSF_{q,m,0}(\la)\cd v\\
	&&=
	\sum_{m=0}^{\iy} x_{\la - \wt(v) -\wt(w)} 
	\ten \FSF_{q,0,m}(\la- \wt(v))\cd w
	\ten \FSF_{q,m,0}(\la)\cd v \\
	&&\qquad+\sum_{m=0}^{\iy} \sum_{m'=1}^{\iy}
	f^{m'}\cd x_{\la - \wt(v) - \wt (w)} \ten 
	\FSF_{q,m',m}(\la - \wt(v))\cd w
	\ten \FSF_{q,m,0}(\la)\cd v\\
&&=\Phi_{q,\la}^u(w\ten v),
\end{eqnarray*}
where
\[
J_{W,V}(\la)(w\ten v)=u=\sum_{m=0}^{\iy}\FSF_{q,0,m}(\la-\wt(v))\cd w
	\ten \FSF_{q,m,0}(\la)\cd v.
\]
Now substitute
\eqref{eq:IM10} and \eqref{eq:IM11}.
$\square$\\

Write the fusion matrix of two Verma modules
$M_{q,\ga}$, $ M_{q,\de}$ as
$J_{q,\de,\ga}(\la):=J_{M_{q,\de},M_{q,\ga}}(\la)$
and write the matrix elements
with respect to their standard bases as
\begin{equation}
	J_{q,\de,\ga} (\la)( x^{\de}_{\de-2s + 2n} 
	\ten 
	x^{\ga}_{\ga-2n})
	= 
	\sum_{m=0}^{n} J_{q,\de,\ga,s;m,n}\;
x^{\de}_{\de- 2s + 2m}\ten x^{\ga}_{\ga - 2m}.
	\label{eq:FM4}
\end{equation}
Then, by comparison with \eqref{eq:FM2} and by change of summation
variable we obtain:
\begin{eqnarray}
 &&\qquad\qquad J_{q,\de,\ga,s;m,n}(\la)
	=
	{\frac{[n]_{q}!}{[m]_{q}!\, [n-m]_{q}!}} 
	q^{-{\frac{n-m}{2}}(\la + 1)}
	q^{{\frac{n-m}{4}}(\de - 2s + \ga)}
	{\frac{([\ga - n + 1 ]_{q})_{n-m} }
	{([-\la + \ga -2n ]_{q})_{n-m}}},\qquad
	\label{eq:FM5}\\
&&\mbox{where}\nonumber\\
&&\qquad\qquad m,n,s\in\ZZ,\quad 0\le m\le n \le s,\quad
\la - \de,
\la - \de- \ga \notin \Zplus.
\label{eq:FM6}
\end{eqnarray}

We will now show that the inverse of the fusion matrix
$J_{q,\de,\ga}$ exists and
that its matrix elements, defined by
\begin{eqnarray}
	J_{q,\de,\ga}^{-1} (\la)( x^{\de}_{\de-2s + 2n} 
	\ten 
	x^{\ga}_{\ga-2n})
	= 
\sum_{m=0}^{n} J_{q,\de,\ga,s;m,n}^{\rm inv}(\la)\;
x^{\de}_{\de- 2s + 2m}\ten x^{\ga}_{\ga - 2m},
	\label{eq:FM7}
\end{eqnarray}
have explicit expression
\begin{equation}
	J_{q,\de,\ga,s;m,n}^{\rm inv}(\la)
	= 
	{\frac{[n]_{q}!}{[m]_{q}!\, [n-m]_{q}!}} 
	q^{-{\frac{n-m}{2}}(\la + 1)}
	 q^{{\frac{n-m}{4}}(\de - 2s + \ga)}
	{\frac{([\ga - n + 1 ]_{q})_{n-m} }
	{([\la - \ga + 2m + 2 ]_{q})_{n-m}}}.
	\label{eq:FM8}
\end{equation}

\noindent
{\bf Proof of (\ref{eq:FM8})}\quad
We have to show that
$\sum_{l=0}^{n-m} J_{q,\de,\ga,s;m,m+l}^{\rm inv}J_{q,\de,\ga,s;m+l,n}(\la)
=\de_{m,n}$ after substitution of
\eqref{eq:FM5} and \eqref{eq:FM8}.
Indeed, by use of the $q$-Chu-Vandermonde 
sum \eqref{eq:PR27} we have:
\begin{eqnarray*}
	\sum_{l=0}^{n-m} J_{q,\de,\ga,s;m,m+l}^{\rm inv}(\la)
J_{q,\de,\ga,s;m+l,n}(\la)
\hskip -0.5truecm
	&&=
	\frac{[n]_{q}!}{[m]_{q}!\,[n-m]_{q}!} q^{-\frac{n-m}{2}(\la +1)}
	q^{\frac{n-m}{4}(\de + \ga - 2s)} 
	\frac{([\ga- n +1]_{q})_{n-m}}
	{([\ga- \la - 2n]_{q})_{n-m}}\\
	&&\qquad\times	
\qhyp21{q^{-n+m},q^{\la -\ga +n +m +1}}{q^{\la - \ga + 2m +2}}{q,q}\\
&&=\frac{[n]_{q}!}{[m]_{q}!\, [n-m]_{q}!} 
	q^{-\frac{n-m}{2}(\la +1)}
	q^{\frac{n-m}{4}(\de + \ga - 2s)}
 q^{(\la-\ga+n+m+1)(n-m)}\\
&&\qquad\times
	 \frac{([\ga- n +1]_{q})_{n-m}}
	{([\ga- \la - 2n]_{q})_{n-m}}\,
	\frac{(q^{m-n+1};q)_{n-m}} 
	{(q^{\la - \ga + 2m +2};q)_{n-m}},
\end{eqnarray*}  
which equals $\de_{m,n}$ because of the factor $(q^{m-n+1};q)_{m-n}$.
$\square$
\section{The universal fusion matrix} \label{sec:UF}
Formula \eqref{eq:FM2} for the fusion matrix
suggests the definition of the {\em  universal fusion matrix},
see \cite{ArBRR},
as a generalized element of $\Uq\ten\Uq$ given by
\begin{equation}
J_q(\la):=
\sum_{l= 0}^{\iy } 
\frac{q^{-{\frac{l}{2}}(\la+1)}}{[l]_{q}!}
f^{l} q^{lh/4}
\ten
e^{l}\frac{q^{lh/4}}{([-\la + h]_{q})_{l}}\,.
\label{eq:UF1}
\end{equation}
This has the property that
\begin{equation}
	J_{q}(\la) \cd (w\ten v)=J_{W,V}(\la) (w\ten v)
\label{eq:UF2}
\end{equation}
for each pair of $\Uq$-modules $W,V$ and for any $w\in W$, $v\in V$.
In fact, $J_q(\la)$ is the unique generalized element in $\Uq\ten\Uq$
of the form
\begin{equation}
J_q(\la)=\sum_{l=0}^\iy J_{q}^{(l)}(\la)
\label{eq:UF6}
\end{equation}
with
$J_{q}^{(0)}(\la)=1\ten 1$ and
$J_{q}^{(l)}(\la)=f^l\phi_l(q^{\quart h},\la)
\ten e^l\psi_l(q^{\quart h},\la)$
($\phi_l$ and $\psi_l$ being rational functions) such that
\eqref{eq:UF2} holds for all pairs $W,V$ of irreducible finite
dimensional $\Uq$-modules.
\begin{definition}
Let $\FSM(\la)$ be a generalized element of 
$\Uq\ten\Uq$ depending on $\la\in\CC$ (outside some discrete subset of $\CC$).
Let $V_1,\ldots,V_n$ be $U_q$-modules such that the action of $\FSM$ on
$V_1\ten\cdots\ten V_n$ is well-defined.
Let $v_j\in V_j$ ($j=1,\ldots,n$) be weight vectors.
Then
\[
\FSM(\la - h^{(i)})\cd
(v_{1} \ten \cdots  \ten v_{n}) :=  
\FSM(\la - \wt(v_{i}))\cd
(v_{1} \ten \cdots \ten v_{n})\qquad
(i\in \{1,\cdots, n\}).
\]
\vspace{-\bigskipamount}
\vspace{-\bigskipamount}
\label{th:UF7}
\end{definition}
\begin{theorem}
The  universal fusion matrix satisfies the identity
\begin{equation}
	\De \FSF_{q,m,n}(\la) J_{q}(\la) =
\sum_{l=0}^\iy \FSF_{q,m,l}(\la - h^{(2)}) \ten \FSF_{q,l,n} (\la),
\label{eq:UF8}
\end{equation}
where $\FSF_{q,m,n}(\la)$ is given by \eqref{eq:IM1}.
\label{th:UF9}
\end{theorem}
\Proof
Let both sides of \eqref{eq:FM1} act on $f^n\cd x_\la$
and use \eqref{eq:IM6} repeatedly. This yields:
\begin{eqnarray*}
&&\sum_{m=0}^\iy f^m\cd x_{\la-\wt(v)-\wt(w)}\ten
\FSF_{q,m,n}(\la)\cd J_{q,W,V}(\la)\,(w\ten v)\\
&&\qquad\qquad=
\sum_{m=0}^\iy f^m\cd x_{\la-\wt(v)-\wt(w)}\ten\sum_{l=0}^\iy
\FSF_{q,m,l}(\la-\wt(v))\cd w\ten\FSF_{q,l,n}(\la)\cd v.
\end{eqnarray*}
By linear independence of the vectors $f^m\cd x_{\la-\wt(v)-\wt(w)}$
\quad($m=0,1,\ldots$) we obtain that
\[
\FSF_{q,m,n}(\la)\cd J_{q,W,V}(\la)\,(w\ten v)
=\sum_{l=0}^\iy
\FSF_{q,m,l}(\la-\wt(v))\cd w\ten\FSF_{q,l,n}(\la)\cd v.
\]
This last identity can be interpreted as identity \eqref{eq:UF8}
acting on $w\ten v$. $\square$\\

Formula \eqref{eq:UF8} implies that the $J_{q}(\la)$ satisfies a
{\em shifted 2-cocycle  condition}
\begin{eqnarray}
	(\id \ten \De)J_{q}(\la)\, (\id \ten J_{q}(\la))
	=
	(\De \ten \id) J_{q}(\la)\, (J_{q}(\la - h^{(3)}) \ten \id).
\label{eq:UF10}
\end{eqnarray}
\begin{remark} \label{th:UF5}
Formula \eqref{eq:UF1} is equivalent to
\cite[formula (3)]{BabBB}. Indeed, the paper \cite{BabBB}
works with generators
$H$, $E_\pm$ satisfying relations
\[
[H,E_\pm]=\pm 2E_\pm,\quad
[E_+,E_-]=\frac{q^H-q^{-H}}{q-q^{-1}}\,,
\]
and comultiplication
\[
\De(H)=H\ten\id+\id\ten H,\quad
\De(E_\pm)=E_\pm\ten q^{\half H}+q^{-\half H}\ten E_\pm\,,
\]
and we can rewrite \cite[(3)]{BabBB} as
\begin{eqnarray}
&&F_{12}(x)=\sum_{k=0}^\iy \frac{x^k}{[k]_{q^2}!}\,q^{\half k H} E_+^k
\ten \frac{(-1)^k q^{\half k H}}{([{}^q{\rm log}\,x+h+k]_{q^2})_k}\,E_-^k
\nonumber\\
&&\qquad\quad=\sum_{k=0}^\iy \frac{x^k}{[k]_{q^2}!}\,E_+^k q^{\half k H}
\ten E_-^k\,\frac{ q^{\half kH}}{([-{}^q{\rm log}\,x-h+1]_{q^2})_k}\,.
\label{eq:UF3}
\end{eqnarray}
The last infinite sum becomes equal to the \RHS\ of
\eqref{eq:UF1} after the successive substitutions
\begin{equation}
q\to q^{-\half},\;
E_+\to \frac{q^{5/4}}{1-q}\,f,\;
E_-\to \frac{1-q}{q^{5/4}}\,e,\;
H\to -h,\;
x\to q^{-\half(\la+1)}.
\label{eq:UF4}
\end{equation}
These substitutions also send the above relations and comultiplication
to the corresponding formulas \eqref{eq:PR33} and \eqref{eq:PR4}.
Furthermore note that these substitutions send the shifted cocycle
condition \cite[(4)]{BabBB} to our formula \eqref{eq:UF10}.
The factor $\frac{q^{5/4}}{1-q}$ and its inverse which appear in  
\eqref{eq:UF4} will play later a role in sending the shifted boundary
in \cite{BabBB} to our formula (see section \ref{sec:SB}).
\end{remark}

We will now give another proof of \eqref{eq:UF10}
by obtaining it as the special case $m=n=0$
of the more general identity
\begin{eqnarray}
	&&(\id \ten \De)(\De \FSF_{q,m,n}(\la) J_{q}(\la)) 
	(\id \ten J_{q}(\la)) \nonumber\\
	&&\qquad\qquad=
	(\De \ten \id) ( \De \FSF_{q,m,n}(\la) J_{q}(\la))
	(J_{q}(\la - h^{(3)}) \ten \id).
\label{eq:UF11}
\end{eqnarray}
{\bf Proof of (\ref{eq:UF11})}\quad
It is sufficient to prove \eqref{eq:UF11} with both sides
acting on any $u\ten v\ten w$, where $u$, $v$ and $w$ are weight vectors
in finite dimensional irreducible $\Uq$-modules $U,V,W$. Then
\eqref{eq:UF11} takes the form
\begin{eqnarray}
	&&(\id \ten \De)(\De \FSF_{q,m,n}(\la) J_{q}(\la)) 
	(\id \ten J_{q}(\la))\cd (u \ten v \ten w) \nonumber\\
	&&\qquad\qquad=
	(\De \ten \id) ( \De \FSF_{q,m,n}(\la) J_{q}(\la))
	(J_{q}(\la - h^{(3)}) \ten \id) \cd (u \ten v \ten w).
\label{eq:UF12}
\end{eqnarray}
For the proof of \eqref{eq:UF12} we will rewrite
both of its sides into expressions \eqref{eq:UF13}
and \eqref{eq:UF14}, respectively, which are equal.
Here we will use \eqref{eq:UF8}
repeatedly.
First the \LHS\ of \eqref{eq:UF12}:
\begin{eqnarray}
	&&(\id \ten \De)(\De \FSF_{q,m,n}(\la) J_{q}(\la))
	\left(u \ten J_{q}(\la) \cd (v \ten w)\right) \nonumber\\
	&&\qquad= 
	(\id \ten \De)( \sum_{t=0}^{\iy}
	 \FSF_{q,m,t}(\la - \wt(J_{q}(\la)(v\ten w))
	\ten \FSF_{q,t,n}(\la)) \left(u \ten J_{q}(\la)\cd 
	(v \ten w)\right) \nonumber\\
	&&\qquad=
	\sum_{t=0}^{\iy} \FSF_{q,m,t}
\left(\la - \wt(v) - \wt(w)\right) \cd u \ten
	\De \FSF_{q,t,n}(\la) J_{q}(\la) \cd (v \ten w)\nonumber\\
	&&\qquad=
	\sum_{t=0}^{\iy} \sum_{k=0}^{\iy}
	\FSF_{q,m,t} \left(\la - \wt(v) - \wt(w)\right) \cd u
	\ten \FSF_{q,t,k} \left(\la - \wt(w) \right)\cd v \ten 
	\FSF_{q,k,n}(\la) \cd w.\qquad
\label{eq:UF13}
\end{eqnarray}
The weight preserving property of 
$J_{q}(\la)$ was used in the second equality.
Next the \RHS\ of \eqref{eq:UF12}:
\begin{eqnarray}
	&&(\De \ten \id)(\De \FSF_{q,m,n}(\la)  J_{q}(\la))
	\left(J_{q}(\la - \wt(w))\cd (u \ten v)\ten w\right) \nonumber\\
	&&\qquad=
	(\De \ten \id)\left(\sum_{k=0}^{\iy} \FSF_{q,m,k}
	(\la - \wt(w)) \ten \FSF_{q,k,n}(\la) \right)	
	\left(J_{q}(\la - \wt(w)) \cd (u \ten v) \ten w \right)\nonumber\\
	&&\qquad=
	\sum_{k=0}^{\iy} \De \FSF_{q,m,k}\left(\la - \wt(w)\right)
	J_{q} \left(\la - \wt(w)\right)\cd (u \ten v) \ten 
	\FSF_{q,k,n}(\la)\cd w\nonumber\\
	&&\qquad=
	\sum_{t=0}^{\iy} \sum_{k=0}^{\iy}
	\FSF_{q,m,t}\left(\la - \wt(w) - \wt(v)\right) \cd u
	\ten \FSF_{q,t,k}\left(\la - \wt(w)\right)\cd v \ten 
	\FSF_{q,k,n}(\la) \cd w.\qquad
\label{eq:UF14}
\end{eqnarray}
Indeed, \eqref{eq:UF13} and \eqref{eq:UF14}
are equal. $\square$
\section{The universal fusion matrix and the shifted boundary}
\label{sec:SB}
Babelon, Bernard \& Billey \cite{BabBB} (see notation of
\cite{BabBB} summarized in Remark \ref{th:UF5})
associate to their universal fusion matrix $F_{12}(x)$
(see \eqref{eq:UF3}) a generalized element $M(x)$ in $\Uq$,
called the {\em shifted boundary} and given by
\begin{equation}
M(x)=\sum_{m,n=0}^\iy
\frac{(-1)^mx^m q^{\half n(n-1)+m(n-m)}}
{[n]_{q^2}!\,[m]_{q^2}!\,\prod_{j=1}^n(xq^j-x^{-1}q^{-j})}\,
E_+^n\,E_-^m\,q^{\half(n+m)H},
\label{eq:SB4}
\end{equation}
such that 
\begin{equation}
	F_{12}(x) =\De M(x)\, \left(\id \ten M(x)\right)^{-1}\,
	\left(M(xq^{H_2}) \ten \id\right)^{-1}.
\label{eq:SB5}
\end{equation}
\par
In the following
we will independently derive an explicit expression for the inverse
of the shifted boundary associated to our 
universal fusion matrix $J_{q}(\la)$ given by \eqref{eq:UF1}.
\begin{theorem}
The universal fusion matrix $J_{q}(\la)$ verifies
\begin{equation}
	\De\left(\FSM_{q}(\la)\right)\,J_{q}(\la)=
	\FSM_{q}(\la - h^{(2)}) \ten \FSM_{q}(\la),
\label{eq:SB1}
\end{equation}
where $\FSM_{q}(\la)$, the inverse shifted boundary, is given by 
\begin{eqnarray}
&&\FSM_{q}(\la) = \sum_{m,n=0}^\iy
\frac{q^{\half mn-\quart m^2}q^{-\half m+\quart n} q^{-\half\la n}}
{\qfac n\, \qfac m\, (1-q)^n}\, f^n e^m\,
\frac{q^{\quart(n-m)h}}{\qpoch{-\la+h}m}
\label{eq:SB2}\\
&&\qquad\qquad
=E_q(q^{\quart-\half\la}fq^{\quart h}) \; \FSA_q(q^{-\la+h},(1-q)^2
 q^{-\frac54-\half\la}eq^{\quart h}),
\label{eq:SB3}
\end{eqnarray}
with $E_q$ and 
$\FSA_q$ respectively given by \eqref{eq:PR36} and \eqref{eq:PR38},
and with the two arguments of $\FSA_q$ satisfying the relation in
\eqref{eq:PR38}.
\end{theorem}
\Proof
Let us first define the element $\wtilde{\FSM}_{q,m,n}(\la) \in \Uq$
 which verifies:
\begin{equation}
	\Phi_{q, \la}^{v}(q^{-n^{2}/4} q^{-nh/4} f^{n}\cd x_{\la})
	=\sum_{m=0}^{\iy}
	 q^{-m^{2}/4} q^{-mh/4} f^{m}\cd x_{\la- wt(v)} \ten
	\wtilde{\FSM}_{q,m,n}(\la)\cd v
\label{eq:SB14}
\end{equation}
By use of the intertwining property of $\Phi_{q, \la}^{v}$ together with
\eqref{eq:IM6} and the property
$h\,\FSF_{q,m,n}(\la)=\FSF_{q,m,n}(\la)\,(h+2m-2n)$ implied by 
\eqref{eq:IM1},
we obtain:
\[
	\wtilde{\FSM}_{q,m,n}(\la) = q^{(n^{2}-m^{2})/4} q^{\la(m-n)/4}
	\FSF_{q,m,n}(\la) q^{-mh/4}.
\]
After substitution of \eqref{eq:IM1} this becomes:
\begin{equation}
	\wtilde{\FSM}_{q,m,n}(\la) = \sum_{j=0}^{m\wedge n}
	q^{(m-j)(n-1)/2} q^{(n^{2}-m^{2})/4} 
	\frac{[n]_{q}!}{[j]_{q}!}\,\frac{(q^{-\la/2}f)^{n-j}
	e^{m-j} q^{(n-j)h/4} q^{-(m-j)h/4}}
	{[n-j]_{q}! [m-j]_{q}! ([-\la + h]_{q})_{m-j}}\,.
\label{eq:SB15}
\end{equation}
It can be shown similarly to the proof of
Theorem \ref{th:UF9} that formula
\eqref{eq:UF8} remains valid if we replace $\FSF$ by
$\wtilde\FSM$. (Just start the proof now by letting both sides of
\eqref{eq:FM1} act on
$q^{-n^{2}/4} q^{-nh/4} f^{n}\cd x_{\la}$ and by applying
\eqref{eq:SB14} repeatedly.)$\;$
Thus we have:
\begin{equation}
	\De\wtilde{\FSM}_{q,m,n}(\la)\,J_{q}(\la) =
	\sum_{l=0}^{\iy}\wtilde{\FSM}_{q,m,l}(\la - h^{(2)}) \ten
	\wtilde{\FSM}_{q,l,n}(\la).
\label{eq:SB16}
\end{equation}
Now put $n:=N-m$ in \eqref{eq:SB16}, substitute
\eqref{eq:SB15}, and sum over $m$ from 0 to $N$:
\begin{eqnarray}
&&\hspace{-3\bigskipamount}\De\Biggl(\sum_{m=0}^N
\sum_{j=0}^{\;m\wedge(N-m)} q^{(m-j)(N-m-1)/2} 
q^{((N-m)^2-m^2)/4}\frac{[N-m]_{q}!}{[j]_{q}!\,[N-m-j]_{q}!\,[m-j]_{q}!}
\nonumber\\
&&\qquad\qquad\qquad\times
	\frac{(q^{-\la/2}f)^{N-m-j} e^{m-j} q^{(N-m-j)h/4} q^{-(m-j)h/4}}
	{([-\la +h]_{q})_{m-j}}\Biggr)J_{q}(\la)
\nonumber\\
&&\hspace{-3\bigskipamount}=
\sum_{m=0}^N	\sum_{l=0}^{\iy}
	\sum_{i=0}^{m \wedge l}
	\sum_{j=0}^{\;l \wedge(N-m)} 
	q^{(m-i)(l-1)/2} q^{(l^2-m^2)/4}\frac{[l]_{q}!}{[i]_{q}!} 
	\frac{(q^{-(\la-h^{(2)})/2}f)^{l-i} e^{m-i} q^{(l-i)h/4} q^{-(m-i)h/4}}
	{[l-i]_{q}! [m-i]_{q}! ([-(\la-h^{(2)}) +h]_{q})_{m-i}}
\nonumber\\	
&&\hspace{-2\bigskipamount}\ten
	q^{(l-j)(N-m-1)/2} q^{((N-m)^2-l^2)/4}\frac{[N-m]_{q}!}{[j]_{q}!} 
	\frac{(q^{-\la/2}f)^{N-m-j} e^{l-j} q^{(N-m-j)h/4} q^{-(l-j)h/4}}
	{[N-m-j]_{q}! [l-j]_{q}! ([-\la +h]_{q})_{l-j}}\,.
\label{eq:SB17}
\end{eqnarray}
The double sum over $m,j$ on the \LHS\ can be rewritten as a double sum
over $m':=m-j$ and $n':=N-m-j$ with $m',n'\ge0$, $m'+n'\le N$ and
$N-m'-n'$ even. Write $m',n'$ again as $m,n$.
Then the \LHS\ of \eqref{eq:SB17} becomes:
\begin{eqnarray*}
&&	\De \Biggl(
\sum_{\begin{array}{c}\\[-1.5\bigskipamount]
\scriptstyle m,n \geq 0\\[-\smallskipamount]
\scriptstyle m+n \leq N\end{array}}
	q^{m(n-m)/4} q^{-m/2} q^{nN/4}\,
	\frac{[\thalf(n-m+N)]_{q}!}{[\thalf(-n-m+N)]_{q}!}\,
	\frac{(q^{-\la/2}f)^{n} e^{m} q^{nh/4} q^{-mh/4}}
	{[n]_{q}!\, [m]_{q}!\, ([-\la +h]_{q})_{m}}\Biggr)J_{q}(\la)\\
&&	=
	\De \Biggl(\sum_{\begin{array}{c}\\[-1.5\bigskipamount]
\scriptstyle m,n \geq 0\\[-\smallskipamount]
\scriptstyle m+n \leq N\end{array}}
	\frac{(q;q)_{\half(n-m+N)}}{(q;q)_{\half(-n-m+N)}}
	(1-q)^{-n}q^{mn/2 - m^2/4 -m/2 +n/4}\,
	\frac{(q^{-\la/2}f)^{n} e^{m} q^{nh/4} q^{-mh/4}}
	{[n]_{q}!\, [m]_{q}!\, ([-\la +h]_{q})_{m}}\Biggr)	J_{q}(\la).
\end{eqnarray*}

The quadruple sum over $m,l,i,j$ on the \RHS\ of of \eqref{eq:SB17}
can be rewritten as a quadruple sum over
$m':=m-i$, $l':=l-i$, $s:=l-j$, $t:=N-m-j$ with $m',l',s,t\ge0$ and
$N-l'-m'+s-t$ even.
Write $m',l'$ again as $m,l$.
Then the \RHS\ of \eqref{eq:SB17} becomes:
\begin{eqnarray*}
&&	\sum_{\begin{array}{c}\\[-1.5\bigskipamount]
\scriptstyle m,l,s,t \geq 0\\[-\smallskipamount]
\scriptstyle -l+m+s+t \leq N\\[-\smallskipamount]
\scriptstyle +l+m-s+t\leq N\end{array}}
\frac{[\half(N+l-m+s-t)]_{q}!\,[\half(N+l-m-s+t)]_{q}!}
	{[\half(N-l-m+s-t)]_{q}!\, [\half(N+l-m-s-t)]_{q}!}
\\[-1.5\bigskipamount]
&&\qquad\qquad\qquad\times
	q^{\frac{m}{2}((N+l-m+s-t)/2-1)} q^{\frac{s}{2}((N+l-m-s+t)/2-1)}
	q^{\frac{N}{4}(l-m-s+t)}
\\
&&\qquad\qquad\qquad\times
\frac{(q^{(-\la +h^{(2)})/2}f)^{l} e^{m} q^{lh/4} q^{-mh/4}}
	{[l]_{q}!\, [m]_{q}!\, ([-\la +h^{(2)}+h]_{q})_{m}}\ten
	\frac{(q^{-\la/2}f)^{t} e^{s} q^{th/4} q^{-sh/4}}
	{[t]_{q}!\, [s]_{q}!\, ([-\la +h]_{q})_{s}}
\end{eqnarray*}
\begin{eqnarray*}
&&	=\hspace{-\bigskipamount}
	\sum_{\begin{array}{c}\\[-1.5\bigskipamount]
\scriptstyle m,l,s,t \geq 0\\[-\smallskipamount]
\scriptstyle -l+m+s+t \leq N\\[-\smallskipamount]
\scriptstyle +l+m-s+t\leq N\end{array}}
	\frac{(q;q)_{\half(N+l-m+s-t)}(q;q)_{\half(N+l-m-s+t)}}
	{(q;q)_{\half(N-l-m+s-t)}(q;q)_{\half(N+l-m-s-t)}}\,
	\frac{q^{lm/2-m^2/4-m/2+l/4}}{(1-q)^l}
	\\[-1.5\bigskipamount]
&&\qquad\qquad\qquad\times
\frac{(q^{(-\la +h^{(2)})/2}f)^{l} e^{m} q^{lh/4} q^{-mh/4}}
	{[l]_{q}!\, [m]_{q}!\, ([-\la +h^{(2)}+h]_{q})_{m}}	\ten
	\frac{q^{ts/2 - s^2/4 -s/2 +t/4}}{(1-q)^t}\,
	\frac{(q^{-\la/2}f)^{t} e^{s} q^{th/4} q^{-sh/4}}
	{[t]_{q}!\, [s]_{q}!\, ([-\la +h]_{q})_{s}}\,.
\end{eqnarray*}

Now consider identity \eqref{eq:SB17} (with both sides rewritten as
above) with $N$ replaced by $2N$ and with $N$ replaced by $2N+1$, add
these two identities, and let $N\to\iy$ in the resulting identity.
We obtain:
\begin{eqnarray*}
&&	\De \Biggl(\sum_{m,n=0}^{\iy}
	\frac{q^{mn/2 - m^2/4 -m/2 +n/4}}{(1-q)^n}\,
	\frac{(q^{-\la/2}f)^{n} e^{m} q^{nh/4} q^{-mh/4}}
	{[n]_{q}!\, [m]_{q}!\, ([-\la +h]_{q})_{m}}\Biggr) J_{q}(\la) \\	 
&&\qquad=	\sum_{m,l=0}^{\iy}	
	\frac{q^{lm/2-m^2/4-m/2+l/4}}{(1-q)^l}\,
	\frac{(q^{(-\la +h^{(2)})/2}f)^{l} e^{m} q^{lh/4} q^{-mh/4}}
	{[l]_{q}!\, [m]_{q}!\, ([-\la +h^{(2)}+h]_{q})_{m}}	\\
&&\qquad\qquad\qquad\qquad\qquad\qquad\ten
	\sum_{t,s=0}^{\iy}
	\frac{q^{ts/2 - s^2/4 -s/2 +t/4}}{(1-q)^t}\,
	\frac{(q^{-\la/2}f)^{t} e^{s} q^{th/4} q^{-sh/4}}
	{[t]_{q}!\, [s]_{q}!\, ([-\la +h]_{q})_{s}}\,.
\end{eqnarray*}
This yields \eqref{eq:SB1} with \eqref{eq:SB2} substituted.
$\square$
\begin{proposition} \label{th:SB18}
The inverse of $\FSM_q(\la)$ formally exists. It is given by
\begin{equation}
	\FSM_q^{-1}(\la) =
	 \frac{(q^{-\la -2};q^{-1})_{\iy}}{(q^{-\la +h -2};q^{-1})_{\iy}} 
	\sum_{m,n=0}^{\iy} \frac{(-1)^{m} (1-q)^{m-2n} 
	q^{\half m^2 - \quart n^2-\half mn +m(-\frac{7}{4}-\half \la) +2n}}
	{\qfac{m}\, \qfac{n}\, \qpoch{\la+2}{n}} f^n e^m q^{\quart(n+m)h}.
\label{eq:SB6}
\end{equation}
\end{proposition}
\Proof
We reason as in \cite[\S7]{R2000a}.
By \eqref{eq:PR36}, \eqref{eq:PR37} and \eqref{eq:PR39}
we see that $\FSM_q(\la)$ is invertible with inverse
\begin{eqnarray*}
\FSM_q^{-1}(\la)= 
	\sum_{m=0}^{\iy} \frac{(-1)^m q^{\half m(m-1)}}{(q;q)_m}
	\frac{(1-q)^{2m} q^{-m(\frac54+\half\la)}}{(q^{-m-1}q^{-\la+h};q)_m}
	 \left(eq^{\quart h}\right)^m
	\sum_{n=0}^{\iy}\frac{(-1)^n q^{n(\quart -\half \la)}
	\left(f q^{\quart h}\right)^n} {(q;q)_n}\,.
\end{eqnarray*}
Then, by use of \eqref{eq:PR34}, we obtain
\begin{eqnarray*}
&&\FSM_q^{-1}(\la) =\sum_{m,n=0}^{\iy}
	\frac{(-1)^n q^{-\quart m^2-\half n^2 +
\half mn-\half m +n(\frac34-\half\la)}}
	{(1-q)^n\, \qfac{m}\, \qfac{n}}\,f^n e^m\,\frac{q^{\quart(n-m)h}}
	{\qpoch{\la -h -2m +2n +2}{m}}\\
	&&\qquad\qquad\qquad \times
	\sum_{k=0}^{\iy} \frac{ q^{-\quart k^2}q^{k(-\half n + \quart - \half \la)}}
	{(1-q)^k\, \qfac{k}}\,
	\frac{\qpoch{-h-m+n}{k}}{\qpoch{\la - h -m +2n+2}{k}}\,.	
\end{eqnarray*}
The inner sum can be formally written as
\begin{eqnarray*}
&&\qhyp11{q^{h+m-n}}{q^{-\la+h+m-2n-2}}{q^{-1},q^{-n-\la-2}}
	= \frac{(q^{-n-\la-2};q^{-1})_{\iy}}{(q^{-\la + h+m-2n-2};q^{-1})_{\iy}}\\
&&\qquad\qquad\qquad\qquad\qquad\qquad\qquad\qquad
	= \frac{(q^{-\la-2};q^{-1})_{\iy} (q^{-\la +h +2m - 2n -2};q^{-1})_{m}}
	{(q^{-\la +h +2m -2n -2};q^{-1})_{\iy} (q^{-\la-2};q^{-1})_{n}},
\end{eqnarray*}
where we used \eqref{eq:PR35}.
This leads to \eqref{eq:SB6} by use of \eqref{eq:PR18}.
$\square$
\begin{lemma} \label{th:SB7}
If $\mu(\la)$ is any solution (like $\FSM_q(\la)$) of \eqref{eq:SB1}
and if $\ga(\la)=\phi(\la-h)/\phi(\la)$ is a formal element not depending on
$e$ and $f$, defined in terms of some nonzero function
$\phi$ of one
variable, then $\mu(\la)\ga(\la)$
satisfies \eqref{eq:SB1}.
\end{lemma}
\Proof
$\De(\ga(\la))$ commutes with $J_q(\la)$ since $h\ten 1+1\ten h$ commutes
with $f^l\ten e^l$ for all $l$. Hence
\begin{eqnarray}
&&\De\bigl(\mu(\la)\,\ga(\la)\bigr)\,J_q(\la)=\De(\mu(\la))\,J_q(\la)\,
\frac{\phi(\la-h\ten1-1\ten h)}{\phi(\la)}\nonumber\\
&&=\Bigl(\mu(\la-h^{(2)})\ten\mu(\la)\Bigr)\,
\frac{\phi(\la-h\ten1-1\ten h)}{\phi(\la-1\ten h)}\,
\frac{\phi(\la-1\ten h)}{\phi(\la)}\nonumber\\
&&=\Bigl(\mu(\la-h^{(2)})\ten\mu(\la)\Bigr)\,
\Bigl(\frac{\phi(\la-h^{(2)}-h)}{\phi(\la-h^{(2)})}\ten
\frac{\phi(\la-h)}{\phi(\la)}\Bigr)
=\mu(\la-h^{(2)})\ga(\la-h^{(2)})\ten\mu(\la)\ga(\la).\nonumber\\
&&\hskip 16truecm\square\nonumber
\end{eqnarray}
\begin{remark}
The successive substitutions \eqref{eq:UF4} send $M(x)$
(given by \eqref{eq:SB4}) exactly to the double sum in \eqref{eq:SB6}.
Hence we have
\begin{equation}
M(x)^{-1}\to\FSM_q(\la)\,
\frac{(q^{-\la-2}; q^{-1})_{\iy}}{{(q^{-\la+h -2};q^{-1})_{\iy}}}
\quad\mbox{(under successive substitutions \eqref{eq:UF4}).}
\label{eq:SB8}
\end{equation}
By Lemma \ref{th:SB7}, the \RHS\ of \eqref{eq:SB8} satisfies
\eqref{eq:SB1} since $\FSM_q(\la)$ does so. This agrees with the
result in \cite{BabBB}
that $M(x)$ satisfies \eqref{eq:SB5}. Indeed, \eqref{eq:SB5} yields
(after the substitutions \eqref{eq:UF4}) equation \eqref{eq:SB1} for
$M(x)^{-1}$.
\end{remark}
\begin{remark}
Rosengren \cite{R2000a}, working with generators $X_+$, $X_-$, $K$, $K^{-1}$
for $\Uq$ which satisfy relations and $*$-structure
\begin{equation}
KX_\pm K^{-1}=q^{\pm\half} X_{\pm},\quad
[X_+,X_-]=\frac{K^2-K^{-2}}{q^\half-q^{-\half}}\,,\quad
K^*=K,\quad
(X_{\pm})^*=-X_\mp\,,
\label{eq:SB12}
\end{equation}
introduced (see \cite[(3.2)]{R2000a}) a generalized element in $\Uq$
given by
\begin{equation}
U_{\la\mu}:=E_q(\mu q^{-\quart}(1-q)X_+ K^{-1})\,
\FSA_q(q\la\mu K^{-4},q^\quart(1-q)\la X_- K^{-1})\,
\frac{(q\la\mu K^{-4};q)_\iy}{(q\la\mu;q)_\iy}.
\label{eq:SB13}
\end{equation}
In his lecture \cite{R2000b} Rosengren next observed
a relationship between Babelon, Bernard \& Billey's shifted boundary 
$M(x)$ and his $U_{\la\mu}$ in $\Uq$, but he did not give the exact
correspondence. In fact, this correspondence is
\begin{equation}
U_{\bar\la,\bar\mu}^*\to M(x)
\quad\mbox{under the successive substitutions given by}
\label{eq:SB9}
\end{equation}
\begin{equation}
q\to q^{2},\; X_-\to E_+, \;
X_+ \to E_-,\; K^{-1}\to q^{\half H},\;
\la\to xq^{-\half},\;\mu\to xq^\half.
\label{eq:SB10}
\end{equation}

If we combine identity \eqref{eq:SB9}
and substitutions \eqref{eq:SB10} 
with identity \eqref{eq:SB8} and substitutions \eqref{eq:UF4} then we
obtain a relationship between Rosengen's generalized element and our
$\FSM_q(\la)$:
\begin{equation}
(U_{\bar\la,\bar\mu}^*)^{-1}\to\FSM_q(\la)\,
\frac{(q^{-\la-2}; q^{-1})_{\iy}}{{(q^{-\la+h -2};q^{-1})_{\iy}}}
\quad\hbox{under the successive substitutions given by}
\label{eq:SB11}
\end{equation}
\begin{equation}
q\to q^{-1},\;
X_-\to\frac{q^{\frac54}}{1-q}f,\;
X_+\to\frac{1-q}{q^{\frac54}}e,\;
K^{-1}\to q^{\quart h},\;
\la\to q^{-\half\la-\quart},\;
\mu\to q^{-\half\la-\frac34}.
\label{eq:SB19}
\end{equation}
This can also be obtained by comparing
\cite[(3.4)]{R2000a} directly with \eqref{eq:SB2}.
\end{remark}
\begin{remark}
Rosengren \cite[(4.8)]{R2000a} gives a generalized conjugation in  $\Uq$
using his element $U_{\la\mu}$. This can be translated by
\eqref{eq:SB11} and \eqref{eq:SB19} into a generalized conjugation
using $\FSM_q(\la)$:
\begin{eqnarray}
&&\bigl(\FSM_q(\la)\bigr)^{-1}
q^{-\quart h} \Biggl(q^{-\half(\la +3)+\quart}
	(1-q)e -
	(1-q)^{-1} q^{-\half (\la-1) - \quart} f\nonumber\\
&&\qquad\qquad+ (q^{-\la -1}+1)\,
	\frac{q^{-\quart h}-q^{\quart h}}{q^{\half} - q^{-\half}} \Biggr) 
	{\cal M}_q(\la)
=	\frac{q^{-\la -1}( q^{\half h}-1) + (q^{-\half h}-1)}
	{q^{\half} - q^{-\half}}\,. \quad
\label{eq:SB20}	
\end{eqnarray}
Note that \eqref{eq:SB20}, in the form with both sides
multiplied on the left by $\FSM_q(\la)$, can also be proved in
a more straightforward way by brute force:
use relations \eqref{eq:PR3}, \eqref{eq:PR18} in order to pull
$q^{-\quart h} e$, $q^{-\quart h} f$ and $q^{-\half h}$ through each term
of the double series \eqref{eq:SB2}.
\end{remark}
\begin{remark}
It is not possible to take limits for $q\to 1$ in formula \eqref{eq:SB2}
for the inverse shifted boundary as it is given there, although
straightforward limit cases for $q\to 1$ are possible for all other
formulas defining or involving the (universal) fusion matrix.
The obstruction for taking limits in \eqref{eq:SB2} is by the factor
$(1-q)^{-n}$.

After rescaling by putting $q^{-\half\la}=c(1-q)$
($c$ constant),
i.e.\ $\la=-2\log(c(1-q))/\log q$,
a limit for $q\to 1$ in \eqref{eq:SB2} becomes possible.
The limit of $\FSM_q(-2\log(c(1-q))/\log q)$ is $\exp(cf)$,
a group element of $SL(2)$. With the same substitution of $\la$,
the universal fusion matrix $J_q(\la)$ given by \eqref{eq:UF1}
tends to $1\ten 1$ as $q\to 1$. Then \eqref{eq:SB1} degenerates
to $\De(\exp(cf))=\exp(cf)\ten\exp(cf)$ (i.e., $\exp(cf)$ is group-like)
and \eqref{eq:SB20} to
\[
\exp(-cf)\,(cf-\thalf h)\exp(cf)=-\thalf h,\quad{\rm i.e.}\quad
\exp({\rm ad}\,cf)(h)=h-2cf.
\]
\end{remark}
\section{The universal fusion matrix and the ABRR equation} \label{sec:AB}
%
Etingof and Schiffmann \cite[Theorem 8.1 and Appendix B]{ES} showed
that the universal fusion
matrix $J_{q}(\la)$  is 
the unique solution of the form \eqref{eq:UF6}
of the equation that they have called the ABRR equation (in reference
to \cite{ArBRR}). They showed that this is also the case when $q=1$,
and used this to compute $J(\la)$ for $\FSU(sl(2))$ (see their example 
after Theorem 8.1). 
Their result coincides with our expression of the universal fusion matrix 
in the classical limit.

In the following we will show directly that our explicit expression of
the universal fusion 
matrix for $\Uq(sl(2))$ verifies the ABRR equation for $\Uq(sl(2))$.

We
first rewrite the ABRR equation for $\Uq(sl(2))$ following our conventions
in section \ref{sec:PR} for the definition of $\Uq(sl(2))$. (This  
slightly differs from the expression given in \cite{ES}, where the
the conventions of \cite[Chapter 6]{CP} are used.)$\;$
Thus the ABRR equation becomes:
\begin{equation}
	J_{q}(\la)  (1 \ten q^{\theta(\la)} ) =
	\FSR_{0}^{21} ( 1\ten q^{\theta(\la)} ) J_{q}(\la)
\label{eq:AB1}
\end{equation}
with $\FSR_0:=\FSR q^{-\quart(h\ten h)}$ and
$\theta(\la):= \thalf(\la + 1)h - \tquart  h^2$.
Since
\[
[h\ten h,q^{\quart jh}e^j\ten q^{-\quart jh}f^j]=
(-2j(h\ten 1)+2j(1\ten h)+4j^2)(q^{\quart jh}e^j\ten q^{-\quart jh}f^j),
\]
we have
\[
q^{\quart(h\ten h)}(q^{\quart jh}e^j\ten q^{-\quart jh}f^j)
q^{-\quart(h\ten h)}=
q^{j^2}\,q^{-\quart jh}e^j\ten q^{\quart jh}f^j,
\]
and thus, by \eqref{eq:PR7},
\begin{equation}
\FSR_0^{21}=\sum_{l=0}^\iy\left(\FSR_0^{(l)}\right)^{21}
\label{eq:AB2}
\end{equation}
with
\begin{equation}
\left(\FSR_0^{(l)}\right)^{21}:=
\frac{(1-q^{-1})^l q^{-\quart l(l-1)} q^{l^2}}
{[l]_q!}\,
q^{\quart lh} f^{l} \ten  q^{-\quart lh}e^{l}.
\label{eq:AB3}
\end{equation}
Substitution of \eqref{eq:UF6} and \eqref{eq:AB2} in
\eqref{eq:AB1} yields:
\begin{eqnarray}
\sum_{n=0}^\iy J_q^{(n)}(\la) &=&
\sum_{l,m=0}^\iy \left(\FSR_0^{(l)}\right)^{21}
(1\ten q^{\half(\la+1)h-\quart h^2}) J_q^{(m)}	
(1 \ten q^{\quart h^2-\half(\la + 1)h}) \nonumber\\
&=&
\sum_{l,m=0}^{\iy} \left(\FSR_0^{(l)}\right)^{21}
(1\ten q^{(\la + 1)m + m^{2} - mh} )
J_{q}^{(m)}(\la).\nonumber
\end{eqnarray}
Hence, by uniqueness of expansion in view of the PBW theorem:
\begin{equation}
J_{q}^{(n)}(\la)=\sum_{l+m=n}
\left(\FSR_0^{(l)}\right)^{21}
(1\ten q^{(\la + 1)m + m^{2} - mh})
J_{q}^{(m)}(\la).
\label{eq:AB4}
\end{equation}
Since $\left(\FSR_0^{(0)}\right)^{21}=1\ten 1$ and
$1\ten q^{(\la + 1)m+m^2-mh}$ is invertible, the terms $J_{q}^{(n)}(\la)$
are uniquely determined by the recurrence \eqref{eq:AB4} together with
the starting value $J_{q}^{(0)}(\la)=1\ten 1$.
In \eqref{eq:UF1} we obtained
\begin{equation}
J_{q}^{(m)}(\la)=
\frac{q^{-\half m(\la+1)}}{[m]_{q}!}\,
f^{m} q^{\quart mh}
\ten
e^{m}\frac{q^{\quart mh}}{([-\la + h]_{q})_{m}}\,.
\label{eq:AB5}
\end{equation}
We will have another proof of \eqref{eq:AB5} if we can show that \eqref{eq:AB4}
is valid after substitution of \eqref{eq:AB3} and \eqref{eq:AB5}.
This is now straightforward.
After the substitutions just mentioned the \RHS\ of \eqref{eq:AB4} becomes:
\begin{eqnarray}
&&\sum_{l+m=n} \frac{(1-q^{-1})^{l}}{[l]_{q}! [m]_{q}!}
q^{-\quart l(l-1)} q^{\half(\la + 1)m - m^2}
q^{-lm} \left(f^{l+m} q^{\quart(l+m)h}
\ten e^{l+m}q^{-\quart(l+3m)h} \frac1{([-\la + h]_{q})_{m}}\right)
\nonumber\\
&&= (1-q^{-1})^{n} q^{-\quart n(n-1)}
\nonumber\\
&&\qquad\times f^{n} q^{\quart nh}\ten e^{n} q^{-\quart nh}
\sum_{m=0}^{n} \frac{(1-q^{-1})^{-m}
q^{\quart m(2n-m-1)} q^{\half(\la+1)m - m^2}
q^{-(n-m)m}}
{[n-m]_{q}! [m]_{q}!
([-\la +h]_{q})_{m}}\,
q^{-\half mh}.\nonumber
\end{eqnarray}
By \eqref{eq:PR2} the last sum equals
\begin{eqnarray}
\frac{1}{[n]_{q}!} \sum_{m=0}^{n} \frac{(q^{-n};q)_{m}}
{(q;q)_{m} (q^{-\la + h};q)_{m}} q^{m} 
&=&
\frac{1}{[n]_{q}!} 
\qhyp21{q^{-n},0}{q^{-\la +h}}{q,q}
\nonumber\\
&=&
\frac{(-1)^{n} q^{n(-\la + h)} 
q^{\half n(n-1)}}{[n]_{q}!\,(q^{-\la + h};q)_{n}}\,,
\nonumber
\end{eqnarray}
where the last equality is obtained by use of the $q$-Vandermonde sum 
\eqref{eq:PR27}.
Hence the \RHS\ of \eqref{eq:AB4} becomes
\[
\frac{(-1)^{n}(1-q^{-1})^n q^{\quart n(n-1)} q^{-n\la}}
{[n]_{q}!}\,
f^{n} q^{\quart nh} \ten 
e^{n} \frac{ q^{\frac34 nh}}{(q^{-\la + h};q)_{n}}\,,
\]
which equals $J_{q}^{(n)}(\la)$ as given by \eqref{eq:AB5}.
\section{The exchange matrix}
\label{sec:EM}
\begin{definition}
Let $V$, $W$ be two $\Uq$-modules, 
$J_{V,W}(\la)$
the fusion matrix, and $\FSR$ the
\mbox{$\FSR$-matrix} \eqref{eq:PR7}.
The {\em exchange matrix} $R_{V,W}(\la)$
is defined by
\begin{equation}
R_{V,W}(\la):= J^{-1}_{V,W}(\la) \FSR^{21} J^{21}_{W,V}(\la)
\label{eq:EM4}
\end{equation}
(see \cite[\S2.2]{ES}), where $J^{21}(\la):= P J(\la) P$ and 
$\FSR^{21} := P \FSR P$.
\label{th:EM5}
\end{definition}
Write the exchange matrix of two Verma modules
$M_{q,\ga}$, $ M_{q,\de}$ as
$R_{q,\ga,\de}(\la):=R_{M_{q,\ga},M_{q,\de}}(\la)$
and write the matrix elements
with respect to their standard bases as
\begin{equation}
R_{q,\ga,\de}(\la) (x^{\ga}_{\ga-2n} \ten x^{\de}_{\de - 2 s + 2n})
=
\sum_{m=0}^{s} R_{q,\ga,\de,s;m,n}(\la)\,
		x^{\ga}_{\ga-2m} \ten
		 x^{\de}_{\de - 2 s + 2m}\,.
\label{eq:EM1}
\end{equation}
Combination of \eqref{eq:EM4} and \eqref{eq:EM1} yields that
\begin{equation}
P\FSR J_{q,\de,\ga}(\la)\,(x_{\de-2s+2n}^\de\ten x_{\ga-2n}^\ga)
=\sum_{m=0}^s R_{q,\ga,\de,s;m,n}(\la)\,J_{q,\ga,\de}(\la)\,
(x_{\ga-2m}^\ga\ten x_{\de-2s+2m}^\de).
\label{eq:EM6}
\end{equation}
From \eqref{eq:EM6} and \eqref{eq:FM6} we see that the following
constraints are required in $R_{q,\ga,\de,s;m,n}(\la)$:
\begin{equation}
m,n,s\in \ZZ,\quad 0\leq m,n \leq s,\quad
\la-\ga,\la -\de,\la -\de-\ga \notin \Zplus.
\label{eq:EM8}
\end{equation}
\begin{theorem}  
The matrix elements of the exchange matrix 
can be expressed in terms of $q$-Racah polynomials as follows:
\begin{eqnarray}
&&R_{q,\ga,\de,s;m,n}(\la)=
 		q^{{\frac{\ga}{4}}(\de - 2s)}
		 q^{\frac{n}{4} (2n + \de + \ga - 2s)}
		q^{\frac{m}{4} (-3\de + \ga +6s - 2m)}
		{\frac{(q^{-\ga};q)_{n}}{
		(q^{\la - \ga + n + 1};q)_{n}}}
		{\frac{(q^{-s}, q^{\de- s +1}; q)_{m}}
		{(q,q^{\de - \la - 2s +m -1};q)_{m}}}\nonumber\\
		&&\qquad\qquad\qquad\qquad\qquad\times 
		R_{m}(\mu(n) , q^{\de- s}, q^{-\la - s -2}, 
		q^{-s-1}, q^{\la - \ga + s+1}).
\label{eq:EM2}
	\end{eqnarray}
Here $\mu(n) := q^{-n} +q^{\la - \ga + n + 1}$, and
	 the $q$-Racah polynomials are given by
\begin{equation}
		R_{m}(\mu(n) , q^{\de- s}, q^{-\la - s -2}, 
		q^{-s-1}, q^{\la - \ga + s+1})
		:=
\qhyp43{q^{-n}, q^{-m}, q^{\la-\ga+n+1},q^{-\la+\de-2s+m-1}}
{q^{-s},	q^{-\ga},q^{\de-s+1}}{q,q}.
\label{eq:EM7}
\end{equation}
\label{th:EM3}
\end{theorem}
\Proof
It follows by successive application of
\eqref{eq:EM6}, \eqref{eq:FM4},
\eqref{eq:PR7}, \eqref{eq:FM7} and
\eqref{eq:FM8} that:
\begin{eqnarray*}
&&\sum_{m=0}^s R_{q,\ga,\de,s;m,n}(\la)\,
(x_{\ga-2m}^\ga\ten x_{\de-2s+2n}^\de)
=\sum_{k=0}^{n} J_{q,\de,\ga,s;k,n}(\la)J^{-1}_{q,\ga,\de}(\la)P\FSR 
	(x^{\de}_{\de -2s + 2k} \ten x^{\ga}_{\ga-2k})
\\
&&=\sum_{k=0}^{n}J_{q,\de,\ga,s;k,n}(\la) J^{-1}_{q, \ga,\de}(\la) P 
	\sum_{j=0}^{\iy} q^{\quart h\ten h} 
	{\frac{(1-q^{-1})^{j}}{[j]_{q}!}} q^{-\frac{j(j-1)}{4}}
(q^{\frac{jh}{4}}e^{j} \ten q^{-\frac{jh}{4}} f^{j})
\cd	(x^{\de}_{\de -2s + 2k} \ten x^{\ga}_{\ga-2k})
\\
&&=\sum_{k=0}^{n}J_{q,\de,\ga,s;k,n}(\la)
\sum_{j=0}^{s-k}\frac{(1-q^{-1})^{j}}{[j]_{q}!}
q^{-\frac{j(j-1)}{4}}
q^{\quart(\de-2s+2k+2j)(\ga-2k-2j)}
q^{\quart j(\de-\ga-2s+4k+4j)}
\\
&&\qquad\qquad\times
(-1)^j([-s+k]_q)_j([\de-s+k+1]_q)_j
J^{-1}_{q, \ga,\de}(\la)\cd (x_{\ga-2k-2j}^\ga\ten x_{\de-2s+2k+2j}^\de)
\\
&&=\sum_{k=0}^{n} \sum_{j=0}^{s-k} \sum_{m=k+j}^{s}
	q^{-\frac{j(j-1)}{4}}
q^{\frac{1}{4}(\de - 2s + 2k+ 2j)(\ga - 2k - 2j)}
		q^{\frac{j}{4}(\de-\ga-2s+4k+4j)}
J_{q,\de,\ga,s;k,n}(\la) J^{\rm inv}_{q,\ga,\de,s;s-m, s-k-j}(\la)
\\
&&\qquad\qquad \times
\frac{(q^{-1}-1)^j([-s+k]_{q})_{j} ([\de - s + k+ 1]_{q})_{j}}{[j]_{q}!}\,
	x^{\ga}_{\ga-2m} \ten x^{\de}_{\de- 2s + 2m}
\end{eqnarray*}
\begin{eqnarray*}
&&=\sum_{k=0}^{n} J_{q,\de,\ga,s;k,n}(\la)
	\sum_{m=k}^{s}
\frac {q^{\quart(\de- 2s+2k)(\ga-2k)} q^{\quart(m-k)(\ga-2s+\de-2\la-2)}
[s-k]_{q}! \, ([\de-s+k+1]_{q})_{m-k}}	
{[s-m]_{q}!\, [m-k]_{q}!\,([\la - \de + 2s - 2 m + 2]_{q})_{m-k}}
\\
&&
\qquad\qquad \times
	 \sum_{j=0}^{m-k} \frac{(-1)^{j}  q^{-\half j(j-1)}
q^{j(-\de-2k+2s+\la+1)}
	(q^{-m + k}, q^{\de - \la - 2s +m +k -1};q)_{j}}
	{(q;q)_{j}}\,
	(x^{\ga}_{\ga - 2m} \ten x^{\de}_{\de - 2s + 2m}).
\end{eqnarray*}
The most inner sum equals
\[
\qhyp20{q^{-m+k}, q^{-\la + \de - 2s + m+k-1}}{-}
{q,q^{-2k - \de + 2s + \la + 1}}
=q^{(k-m)(-\la+\de-2s+m+k-1)}
\]
by the limit case of the $q$-Chu-Vandermonde sum \eqref{eq:PR29}.
Now substitute \eqref{eq:FM5} and interchange the summations
over $k$ and $m$. Then, in view of \eqref{eq:FM5}, we obtain:
\begin{eqnarray*}
	R_{q,\ga,\de,s;m,n}(\la)
	&=&
	\frac{q^{\frac{\ga}{4}(\de - 2 s)}
	q^{\frac{n}{4}(\de + \ga - 2 \la - 2s -2 )}
	q^{\frac{m}{4}( 2 \la - 3 \de + \ga + 6s - 4m +2)}
	([-s]_{q})_{m} ([\de - s +1]_{q})_{m}([-\ga]_{q})_{n}}
	{[m]_{q}!\, ([-\la + \de - 2s + m -1]_{q})_{m}
	([\la - \ga + n +1]_{q})_{n}}
	\\
	&&\qquad\times \sum_{k=0}^{m \wedge n} 
	\frac{([-n]_{q})_{k}([-m]_{q})_{k}
	([\la-\ga+n+1]_{q})_{k} 
	([-\la + \de - 2s +m -1]_{q})_{k}}
	{[k]_{q}!\, ([-s]_{q})_{k} ([-\ga]_{q})_{k} 
	([\de - s + 1]_{q})_{k}}\,,
\end{eqnarray*}
which can be rewritten as
\eqref{eq:EM2}.
$\square$\\

Note that in the above proof, the $j$-sum and its evaluation would not 
occur in the corresponding $q=1$ case (the exchange matrix is defined then
by $R_{V,W}(\la):= J^{-1}_{V,W}(\la) J^{21}_{W,V}(\la)$).
\section{The exchange matrix and q-Racah coefficients}
\label{sec:ER}
The $q$-Racah coefficients arise as in the classical case when one considers
two different ways of decomposing the tensor product of three 
finite dimensional irreducible representations of 
$\Uq$ (see for example \cite{BiLoh} for the $q$-case and 
\cite{BiLou} for the $q=1$ case).
Use the notation of Remarks \ref{th:PR20} and \ref{th:OS2}.
Let $j_1,j_2,j_3,j\in\thalf\Zplus$ be such that
\begin{equation}
j_1+j_2+j_3-j,\;
j_1+j_2-j_3+j,\;j_1-j_2+j_3+j,\;-j_1+j_2+j_3+j\in\Zplus\,,
\label{eq:ER7}
\end{equation}
and let $j_{12}\in\thalf\Zplus$ be such that
\begin{equation}
 |j-j_3| \vee |j_1-j_2|\, \leq \, j_{12} \, \leq \,(j+j_3) \wedge (j_1+j_2).
\label{eq:ER8}
\end{equation}
Inequalities \eqref{eq:ER7} give precisely the condition that
$V^j$ occurs at least once
in $V^{j_1}\ten V^{j_2}\ten V^{j_3}$.
Furthermore, \eqref{eq:ER7} combined with \eqref{eq:ER8} is precisely
the condition that
$V^{j_{12}}$ occurs in $V^{j_1}\ten V^{j_2}$
and that
$V^j$ occurs in $V^{j_{12}}\ten V^{j_3}$.
Let $e_m^{j_{12},j}(j_1,j_2\mid j_3)$ ($m=-j,-j+1,\ldots,j$)
form the standard basis of $V^j$ within $V^{j_{12}}\ten V^{j_3}$ within
$V^{j_1}\ten V^{j_2}\ten V^{j_3}$, where we twice used the normalization of
a standard basis of an irreducible submodule of a tensor product as
described in Remark \ref{th:OS2}.

If we combine this definition of $e_m^{j_{12},j}(j_1,j_2\mid j_3)$ with
formula \eqref{eq:OS7}, applied twice, then we obtain:
\begin{equation}
e_m^{j_{12},j}(j_1,j_2\mid j_3)=
(\wtilde{\Phi}_{q,j_{12}}^{e^{j_2}_{j_{12} - j_1}}
	\ten \id) \circ
	\wtilde{\Phi}_{q,j}^{e^{j_3}_{j-j_{12}}}(e_m^j).
\label{eq:ER3}
\end{equation}

Similarly as above,
let $e_m^{j_{23},j}(j_1\mid j_2,j_3)$ ($m=-j,-j+1,\ldots,j$)
form the standard basis of
$V^j$ within $V^{j_{1}}\ten V^{j_{23}}$ within
$V^{j_1}\ten V^{j_2}\ten V^{j_3}$.
Then the {\em $q$-Racah coefficients}
 ${}_qW^{j_2,j_3,j_{23}}_{j_{12}-j_1,j-j_{12},j-j_1}(j)$
are defined by:
\begin{equation}
e_m^{j_{12},j}(j_1,j_2\mid j_3)=
\sum_{j_{23}}{}_qW^{j_2,j_3,j_{23}}_{j_{12}-j_1,j-j_{12},j-j_1}(j)\,
e_m^{j_{23},j}(j_1\mid j_2,j_3)
\end{equation}
where we sum over
\begin{equation}
|j-j_1| \vee |j_2-j_3|\, \leq \, j_{23} \, \leq \,(j+j_1) \wedge (j_2+j_3).
\label{eq:ER88}
\end{equation}
The notation we use for $q$-Racah coefficients is
defined in \cite[Definition 3.72]{BiLoh}.

Graphically (see \cite{KiR}), the definition of $q$-Racah coefficients
has the form\\

\hspace{2cm} \includegraphics[height=2.5cm]{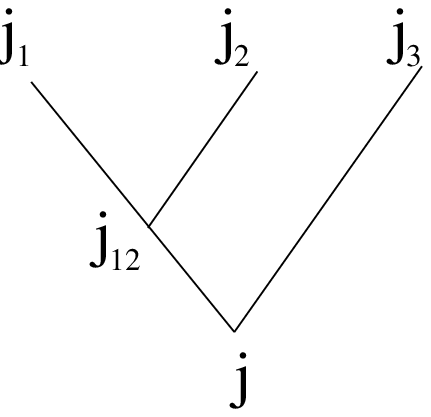}

\vspace{-1.8cm}

\hspace{4.5cm}$= \sum_{j_{23}} 
{}_qW^{j_2,j_3,j_{23}}_{j_{12}-j_1,j-j_{12},j-j_1}(j)$

\vspace{-1.4cm}

\hspace{10cm}\includegraphics[height=2.5cm]{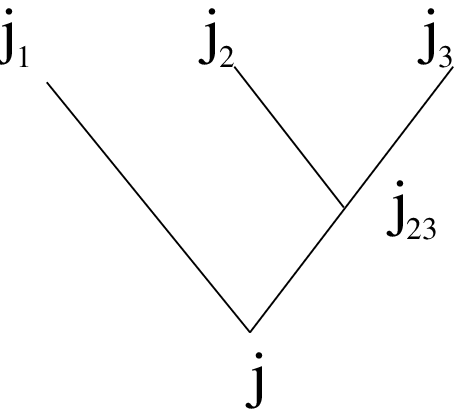}

\vspace{-2.5cm}

\begin{equation}
\label{eq:ER12}
\end{equation}\\
\\
Here, for the diagram both on the left and on the right one has to
substitute the $m$th standard basis vector in the
module isomorphic to $V^j$ which is evidently determined within
$V^{j_1}\ten V^{j_2}\ten V^{j_3}$ by the corresponding diagram.

Because the choice of $m$, above and in the sequel, is irrelevant,
we do not need to put $m$ in the diagram. Therefore, unlike as in
\cite{KiR}, we put $j$-labels
on the vertices of the diagrams and we do not label the edges.
\\[\medskipamount]
\indent
We will need the following formula stated in \cite[(5.11)]{KiR}.
\begin{eqnarray}
&&(P\FSR)_{23}\,e_m^{j_{13},j}(j_3,j_1\mid j_2)=
\sum_{j_{12}} (-1)^{j_{12} + j_{13} -j -j_1}
q^{\half(c_j + c_{j_1} - c_{j_{13}} -c_{j_{12}})}\nonumber\\
&&\qquad\qquad\qquad\qquad\qquad\qquad\times ~
{}_qW^{j_1,j_2,j_{12}}_{j_{13}-j_{3}, j-j_{13},j-j_3}(j)\,
e_m^{j_{12},j}(j_1,j_2\mid j_3).
\label{eq:ER13}
\end{eqnarray}
Here $c_j:=j(j+1)$ as before.
Formula \eqref{eq:ER13} can be written graphically as follows.

\vspace{1cm}

\includegraphics[height=2.5cm]{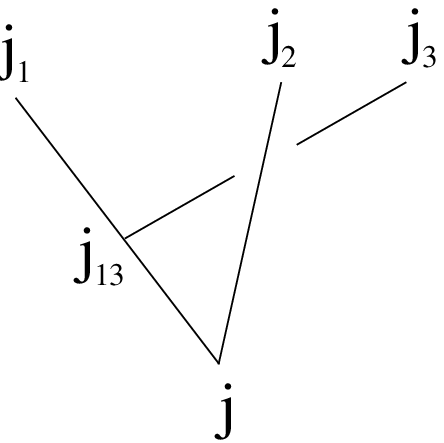}

\vspace{-1.5cm}

\hspace{2cm} $= \sum_{j_{12}} (-1)^{j_{12} + j_{13} -j -j_1}
q^{\half(c_j + c_{j_1} - c_{j_{13}} -c_{j_{12}})}
{}_qW^{j_1,j_2,j_{12}}_{j_{13}-j_{3}, j-j_{13},j-j_3}(j)$

\vspace{-1.6cm}

\hspace{13.2cm}\includegraphics[height=2.5cm]{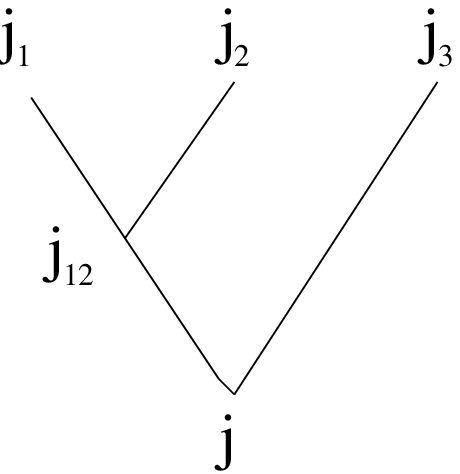}

\vspace{-1.5cm}

\begin{equation}
\label{eq:ER14}
\end{equation}\\
The interpretation of the diagrams is similar as we
explained for \eqref{eq:ER12}.

\smallskip\noindent{\bf Proof of (\ref{eq:ER13})}\quad
Rewrite \eqref{eq:ER12} as

\hspace{2cm}\includegraphics[height=2.5cm]{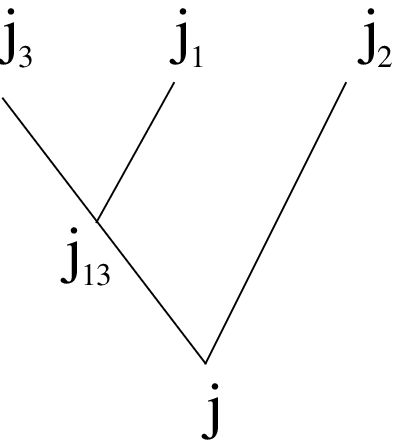}

\vspace{-1.5cm}

\hspace{4cm}$= \sum_{j_{12}} 
{}_qW^{j_1,j_2,j_{12}}_{j_{13}-j_{3}, j-j_{13},j-j_3}(j)$

\vspace{-1.7cm}

\hspace{9cm}\includegraphics[height=2.5cm]{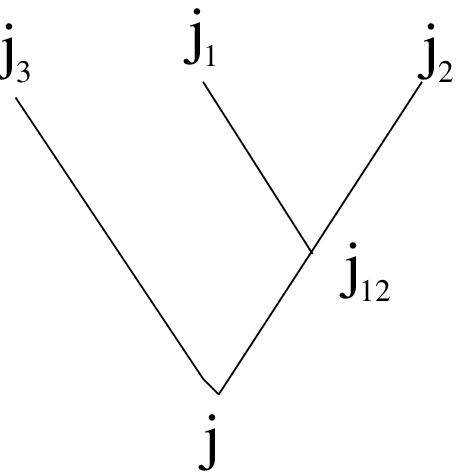}

\vspace{-2.5cm}

\begin{eqnarray}
\label{eq:ER15}
\end{eqnarray}
\\
\\
Then, by the first identity in \eqref{eq:PR22},
we have
\begin{equation}
(P\FSR)_{23}(P\FSR)_{12}=P_{23}P_{12}\FSR_{13}\FSR_{12}=
P_{1,23}(\id\ten\De)\FSR=
(\De\ten\id)(P\FSR).
\label{eq:ER16}
\end{equation}
If we let the first and last part of \eqref{eq:ER16} act on a threefold
tensor product, then this implies an operator identity
\begin{equation}
(P\FSR)_{23}\circ(P\FSR)_{12}=(P\FSR)_{1,23}\,.
\label{eq:ER17}
\end{equation}
If we let the \LHS\ and the \RHS\ of \eqref{eq:ER17} respectively
act on the \LHS\ and \RHS\ of \eqref{eq:ER15} then we obtain by
twofold application of \eqref{eq:AP3}:
\vspace{1.5cm}

\hspace{-.5cm} 
$(-1)^{j_{13}-j_1-j_3} ~ q^{\half(c_{j_{13}}-c_{j_1}-c_{j_3})}~
\left( P \cal{R} \right)_{23}$
\vspace{-1.5cm}

\hspace{6cm}\includegraphics[height=2.5cm]{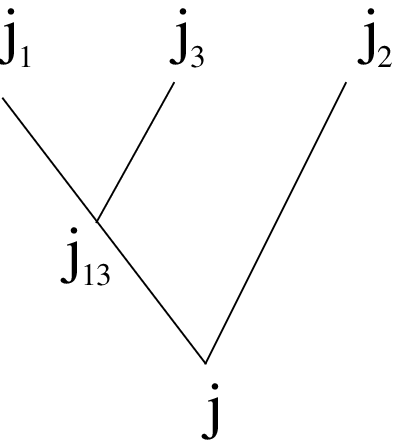}

\vspace{-1.7cm}

\hspace{7.7cm}$=\sum_{j_{12}}~ (-1)^{j-j_3-j_{12}}~
q^{\half(c_j-c_{j_3}-c_{j_{12}})} $

\vspace{-1.5cm}

\hspace{13.7cm}\includegraphics[height=2.5cm]{A5.2.eps}

\vspace{-1cm}
\hspace{8cm} $ \times~{ }_qW^{j_1,j_2,j_{12}}_{j_{13}-j_{3}, j-j_{13},j-j_3}(j)$

\vspace{-.5cm}
\begin{equation}
\label{eq:ER18}
\end{equation}
This is equivalent to \eqref{eq:ER14}.
$\square$

\begin{theorem}
\label{th:ER4}
The $q$-Racah coefficients can be expressed in terms of
matrix elements of the exchange matrix as follows:
\begin{eqnarray}
&&(-1)^{j+j_1-j_{12}-j_{13}}
q^{\half(c_j+c_{j_1}-c_{j_{12}}-c_{j_{13}})}
{}_qW^{j_1,j_2,j_{12}}_{j_{13}-j_{3}, j-j_{13},j-j_3}(j)
=\wtilde{R}^{j_2 j_3}_{q,{\scriptscriptstyle j_{12}-j_1,j-j_{12};
	j-j_{13},j_{13}-j_1}} (j )
\nonumber
\\
&&
:=\frac
{\|\Phi_{q,2j_{12}}^{x_{2j_{12}-2j_1}^{2j_2}}(x_{2j_{12}})\|\,
\|\Phi_{q,2j}^{x_{2j-2j_{12}}^{2j_3}}(x_{2j})\|}
{\|\Phi_{q,2j_{13}}^{x_{2j_{13}-2j_1}^{2j_3}}(x_{2j_3})\|\,
\|\Phi_{q,2j}^{x_{2j-2j_{13}}^{2j_2}}(x_{2j})\|}\,
R_{q,\,2j_2,\,2j_3,\,j_1+j_2+j_3-j;\;j_1+j_2-j_{12},\,j_2+j_{13}-j}(2j).
\qquad\quad
\label{eq:ER5}
\end{eqnarray}
Here $j_1,j_2,j_3,j,j_{12},j_{13}\in\thalf\Zplus$
are constrained  by \eqref{eq:ER7}, \eqref{eq:ER8}
and by \eqref{eq:ER19} below:
\begin{equation}
 |j-j_2| \vee |j_3-j_1|\, \leq \, j_{13} \, \leq \,(j+j_2) \wedge (j_3+j_1).
\label{eq:ER19}
\end{equation}
Furthermore, $c_j:=j(j+1)$.
\end{theorem}
\Proof
First observe from \eqref{eq:EM6}, \eqref{eq:IM28} and
\eqref{eq:FM1}
that
\begin{equation}
	(\id \ten P {\cal R})\circ
	(\Phi_{q,\la - \ga +2n}^{x_{\de-2s+2n}^{\de}} \ten \id)\circ
	\Phi_{q,\la}^{x_{\ga-2n}^\ga}
	=\sum_{m=0}^s
	 R_{q,\ga,\de,s;m,n}(\la)\,
	({\Phi}_{q,\la - \de +2s - 2m}^{ x_{\ga-2m}^\ga}
	\ten\id)\circ
	\Phi_{q,\la}^{x_{\de-2s+2m}^\de}.
\label{eq:ER20}
\end{equation}
Put
\[
\la=2j,\;
\ga=2j_2,\;
\de=2j_3,\;
s=j_1+j_2+j_3-j,\;
n=j_{13}+j_2-j.
\]
It follows by Remark \ref{th:OS2} that formula \eqref{eq:ER20}
remains valid for $j_1,j_2,j_3,j,j_{13}$
being constrained as stated in the theorem.
(Note that the constraints \eqref{eq:ER7}, \eqref{eq:ER8}
and \eqref{eq:ER19} are obtained from the constraints
\eqref{eq:ER7}, \eqref{eq:ER88} and \eqref{eq:ER8}, respectively,
by replacing $j_1,j_2,j_3,j_{12},j_{23}$ by
$j_3,j_1,j_2,j_{13},j_{12}$.)$\;$
Formula \eqref{eq:ER20} then becomes:
\begin{eqnarray*}
	&& (\id \ten P\FSR)\circ
	(\Phi_{q,2j_{13}}^{x_{2j_{13}-2j_1}^{2j_3}}\ten\id)\circ 
	\Phi_{q,2j}^{x_{2j - 2j_{13}}^{2j_2}}
	\\
	&&\quad
	=\sum_{j_{12} = |j-j_3| \vee |j_1 -j_2|}^
	{(j+j_3) \wedge (j_1 +j_2)}
	R_{q,\,2j_2,\,2j_3,\,j_1+j_2+j_3-j\,;\;j_1+j_2-j_{12},\,j_2+j_{13}-j}(2j)\,
	(\Phi_{q,2j_{12}}^{x^{2j_2}_{2j_{12}-2j_1}}\ten\id) \circ
	\Phi_{q,2j}^{x^{2j_3}_{2j-2j_{12}}}.\qquad
\end{eqnarray*}
Next substitute \eqref{eq:OS6} four times in this last identity. Then:
\begin{eqnarray*}
&&(\id \ten P\FSR)\circ
(\wtilde\Phi_{q,j_{13}}^{e_{j_{13}-j_1}^{j_3}}\ten\id)\circ 
\wtilde\Phi_{q,j}^{x_{j -j_{13}}^{j_2}}=
\sum_{j_{12} = |j-j_3| \vee |j_1 -j_2|}^{(j+j_3) \wedge (j_1 +j_2)}
\frac
{\|\Phi_{q,2j_{12}}^{x_{2j_{12}-2j_1}^{2j_2}}(x_{2j_{12}})\|\,
\|\Phi_{q,2j}^{x_{2j-2j_{12}}^{2j_3}}(x_{2j})\|}
{\|\Phi_{q,2j_{13}}^{x_{2j_{13}-2j_1}^{2j_3}}(x_{2j_3})\|\,
\|\Phi_{q,2j}^{x_{2j-2j_{13}}^{2j_2}}(x_{2j})\|}
\qquad
\\
&&\qquad\qquad\qquad\qquad\times
R_{q,\,2j_2,\,2j_3,\,j_1+j_2+j_3-j;\;j_1+j_2-j_{12},\,j_2+j_{13}-j}(2j)\,
(\wtilde\Phi_{q,j_{12}}^{x^{j_2}_{j_{12}-j_1}}\ten\id)\circ
\wtilde\Phi_{q,j}^{x^{j_3}_{j-j_{12}}}.
\end{eqnarray*}
Let both sides of this identity act on $e_m^j$, then
substitute \eqref{eq:ER3} twice in the resulting identity. 
This yields:
\begin{eqnarray*}
&&(P\FSR)_{23}\,e_m^{j_{13},j}(j_3,j_1\mid j_2)=
\sum_{j_{12} = |j-j_3| \vee |j_1 -j_2|}^{(j+j_3) \wedge (j_1 +j_2)}
\frac
{\|\Phi_{q,2j_{12}}^{x_{2j_{12}-2j_1}^{2j_2}}(x_{2j_{12}})\|\,
\|\Phi_{q,2j}^{x_{2j-2j_{12}}^{2j_3}}(x_{2j})\|}
{\|\Phi_{q,2j_{13}}^{x_{2j_{13}-2j_1}^{2j_3}}(x_{2j_3})\|\,
\|\Phi_{q,2j}^{x_{2j-2j_{13}}^{2j_2}}(x_{2j})\|}
\\
&&\qquad\qquad\times
R_{q,\,2j_2,\,2j_3,\,j_1+j_2+j_3-j;\;j_1+j_2-j_{12},\,j_2+j_{13}-j}(2j)\,
e_m^{j_{12},j}(j_1,j_2\mid j_3).
\end{eqnarray*}
Finally compare
with formula \eqref{eq:ER13}. $\square$\\

In the literature and similar to the classical case,
{\em $q$-$6j$ coefficients} are defined in terms of
$q$-Racah coefficients (see either \cite[\S 3.6.1]{BiLoh} or 
\cite[p. 307]{KiR}):
\begin{equation}
\qsixj {j_3} {j_1} {j_{13}} {j_2} {j} {j_{12}}
:=(-1)^{j_1+j_2+j+j_3} \left([2j_{12}+1]_{q} [2j_{13}+1]_{q}\right)^{-\half}
{}_qW^{j_1,j_2,j_{12}}_{j_{13}-j_{3}, j-j_{13},j-j_3}(j).
\label{eq:ER9}
\end{equation}

An explicit expression of $q$-$6j$ coefficients in terms of
${}_4\phi_3$ $q$-hypergeometric functions has been derived in 
\cite{KaK} and in \cite{KiR}. We use here the expression as a finite sum
given in 
\cite[(3.69)]{BiLoh}:
\begin{eqnarray}
&&\qsixj {j_3} {j_1} {j_{13}} {j_2} {j} {j_{12}}
	= \De(j_1,j_2,j_{12})\, \De(j_3,j_1,j_{13})\, \De(j,j_{12}, j_3)
	\, \De(j, j_2,j_{13})
	\nonumber\\[\medskipamount]	
&&\qquad\qquad\times 
	\sum_n \frac{ (-1)^{n} \qfac{n+1}}{\qfac{n-j_1-j_2-j_{12}}
	\qfac{n-j_1-j_3-j_{13}} \qfac{n-j-j_2-j_{13}} 
	\qfac{n-j_3-j-j_{12}}}\nonumber\\
&&\qquad\qquad
	\times \frac{1}{ \qfac{j_1+j_2+j_3+j-n} \qfac{j_2+j_3+j_{12}+j_{13}-n}
	\qfac{j_1+j+j_{12}+j_{13}-n}}
	\label{eq:ER10}
\end{eqnarray}
	where the summation range is
\begin{eqnarray}
\max (j_1 +j_2 +j_{12}, j_1 +j_3 + j_{13},j_2+j+j_{13}, j_3+j+j_{12})
	\leq  n \leq \nonumber \\
	 \min(j_1+j_2+j_3+j, 
	 j_2+j_3+j_{12}+j_{13},j_1+j+j_{12}+j_{13})
\label{eq:ER21}
\end{eqnarray}
and where
\begin{equation}
	\De(j_1,j_2,j):= \left(\frac{ \qfac{-j_1+j_2+j}\, \qfac{j_1-j_2+j}\,
	 \qfac{j_1+j_2-j}}
	{\qfac{j_1+j_2+j+1}} \right)^\half.
\label{eq:ER11}
\end{equation}
Formula \eqref{eq:ER10} may be rewritten as  a
${}_4\phi_3$ $q$-hypergeometric function, depending on certain inequalities
involving the parameters.

\begin{remark}
Equation \eqref{eq:ER5} connects
 $_qW^{j_1,j_2,j_{12}}_{j_{13}-j_{3}, j-j_{13},j-j_3}(j)$
on its \LHS\ with \\
$R_{q,\,2j_2,\,2j_3,\,j_1+j_2+j_3-j;\;j_1+j_2-j_{12},\,j_2+j_{13}-j}(2j)$
on its \RHS. The expression on the \LHS\ can be written as a finite sum
by \eqref{eq:ER9} and \eqref{eq:ER10}, where the summation bounds depend
on certain inequalities involving the parameters.
The expression on the \RHS\ can be written as a limit case of the
finite $q$-hypergeometric sum \eqref{eq:EM7}, where the limit
will also depend on certain inequalities involving the parameters.
The two sums on the two sides of \eqref{eq:ER5} were derived in very
different ways, but must be equal to each other because of the truth of
Theorem \ref{th:ER4}, which we proved in a conceptual way.
Let us independently verify that the two sums are equal.

Let $\ga,\de,\la,s,m,n$ satisfy the constraints \eqref{eq:EM8}.
By use of Sears' transformation \eqref{eq:PR31}, the
$q$-hypergeometric expression \eqref{eq:EM2}
of the exchange matrix can be rewritten as follows
(the apparent singularity for $m>n$ is in fact removed):
\begin{eqnarray*}
&&R_{q,\ga,\de,s;m,n}(\la)=
 	q^{{\frac{\ga}{4}}(\de - 2s)}
q^{\frac{n}{4} (2n + \de + \ga - 2s)} q^{\frac{m}{4} (-3\de + \ga +6s - 2m)}
	{\frac{(q^{-\ga};q)_{n}}{(q^{\la - \ga + n + 1};q)_{n}}}
	{\frac{(q^{-s}, q^{\de- s +1}; q)_{m}}
{(q,q^{\de - \la - 2s +m -1};q)_{m}}}\\
&&\qquad 
	\times	
	\frac{(q^{-n}, q^{s-m+n+\la-\de-\ga+1},q^{-\la - s -1};q)_{m}}
	{(q^{-\ga}, q^{-s}, q^{s-m-\de};q)_{m}}	
	\qhyp43{q^{-m}, q^{n-\ga}, q^{n-s},q^{s-m-\de}}
	{q^{s-m+n+\la-\de-\ga+1},q^{-\la -s -1},q^{n-m+1}}{q,q}
\end{eqnarray*}
\begin{eqnarray*}
&&\qquad = q^{{\frac{\ga}{4}}(\de - 2s)}
	q^{\frac{n}{4} (2n + \de + \ga - 2s)}
	q^{\frac{m}{4} (-3\de + \ga +6s - 2m)} 
	q^{\frac{n}{2}(-\la -n -1)} q^{\frac{m}{2}(\la -m+1)}\\
&&\qquad
	\times \sum_{t=0 \vee (m-n)}^{m \wedge (s-n)}
	\frac{ \qpoch {-\ga} n \qpoch{n-\ga}t}{\qfac{t} \qpoch{-\ga} m}
	\frac{\qpoch{n-m+1}m}{\qpoch{n-m+1}t}
	\frac{\qpoch{s-m+n +\la -\de - \ga +1}m}
	{\qpoch{s-m+n +\la -\de - \ga +1}t}\\
&&\qquad
	\times
	\frac{\qpoch{-\la-s-1}m}{\qpoch{-\la-s-1}t}
	\frac{\qpoch{-m}t \qpoch{n-s}t \qpoch{s-m-\de}t}
	{\qfac{m} \qpoch{\de -\la -2s +m -1}m \qpoch{\la - \ga +n +1}n}
\end{eqnarray*}
If we substitute 
$j_1,j_2,j_3,j,j_{12},j_{13}$ as  in Theorem \ref{th:ER4} and 
take limits to the constraints for $j_1,j_2,j_3,j,j_{12},j_{13}$
as assumed there, and if we 
pass to the new 
summation variable $k:= j+j_1+j_2+j_3-t$, then we obtain
\begin{eqnarray}
&&	R_{q,2j_2,2j_3,j_1+j_2+j_3-j;j_1-j_{12}+j_2,j_{13}-j+j_2} (2j)=
(-1)^{j_{13} +j_{12} +j_2 +j_3} q^{\frac{j_1}{2}(-j_{13} + j_1 +j_{12} +1)}
	q^{\frac{j}{2}(j+j_{12}-j_{13}+1)} \nonumber\\
	[\bigskipamount]
&&\qquad	
	\times q^{-j_{12}(j_{12} +\half)} q^{-\frac{j_{13}}{2}}
	 \frac{\qfac{j_{13} -j_2 +j} \qfac{j_{13} + j_{2}-j}
	\qfac{j_1-j_3 +j_{13}} \qfac{2j_{12} +1}}
	{\qfac{2j_{13}} \qfac{j + j_3 +j_{12} +1}
\qfac{j_{12} +j_1 +j_2 +1}}\nonumber\\
	[\bigskipamount]
&&\qquad
	\times \qfac{j_2 -j_1 +j_{12}} \qfac{j_3 +j -j_{12}}
\qfac{j_1 + j_3 -j_{13}}\nonumber\\
	[\bigskipamount]
&&\qquad
\times \sum_k  
	\frac{(-1)^k \qfac{k+1}}{\qfac{j_{13}+j_{12}+j_2+j_3-k}
	\qfac{j_{13}+j_{12}+j+j_1-k} \qfac{-j-j_3-j_{12}+k}}\nonumber\\
&&\qquad
	\times \frac{1}{ \qfac{-j_1 -j_3 -j_{13} +k}  \qfac{-j_1-j_2-j_{12}+k}
	\qfac{-j_{13}-j_2-j+k} \qfac{j+j_1+j_2+j_3-k}},\nonumber\\
\label{eq:ER22}
\end{eqnarray}
where the summation range is as in \eqref{eq:ER21}.

Now the \LHS\ of \eqref{eq:ER5} with \eqref{eq:ER9} and \eqref{eq:ER10}
substituted equals
the \RHS\ of \eqref{eq:ER5} with \eqref{eq:ER22} and
\eqref{eq:OS5} substituted.
\end{remark}
\section{QDYBE and q-Racah coefficients}
\label{sec:YB}
The exchange matrix satisfies the
quantum dynamical Yang-Baxter equation (QDYBE)
\begin{equation}
{R}_{q}^{23}(\la){ R}_{q}^{13}(\la-h^{(2)})
{R}_{q}^{12}(\la) 
=
{R}_{q}^{12}(\la-h^{(3)}) {R}_{q}^{13}(\la) 
{R}_{q}^{23}(\la-h^{(1)}),
\label{eq:YB1}
\end{equation}
see \cite[Proposition 2.4 and \S2.2]{ES} and \cite[\S4]{Koo2000a}
(which gives some details of proofs in \cite{ES}).
In the following, we will take the limit of the QDYBE \eqref{eq:YB1}
to the case of finite dimensional irreducible representations.
Together with \eqref{eq:ER5} and \eqref{eq:ER9}
this will yield an identity for sums of products of three
$q$-$6j$ symbols, earlier known in the literature and expressing a
symmetry of $q$-$9j$ symbols \cite[(8.1)]{N1989}. 
This consequence of the QDYBE was earlier
mentioned, without giving details, in \cite[\S 8]{EV}.

Let $e_{j_4- j_5}^{j_1}$,$ e_{j_5- j_6}^{j_2}$, 
$e_{j_6 - j_7}^{j_3}$ be  basis vectors 
of the finite dimensional irreducible $\Uq$-modules 
$V^{j_1}$, $V^{j_2}$, $V^{j_3}$. 
The action of a limit case of the QDYBE on $(e_{j_4- j_5}^{j_1} \ten
 e_{j_5- j_6}^{j_2} \ten e_{j_6 - j_7}^{j_3})$
is then given by the following identity:
\begin{eqnarray}
	&&\sum_{j_8,j_9,j_{10}}
	{R}_{q,2j_2, 2j_3, j_2 +j_3 +j_8 -j_4;j_2+j_8-j_9,j_2+j_{10}-j_4}
	(2j_4)
	 {R}_{q,2j_1,2 j_3, j_1 +j_3 +j_7 -j_{10};j_1 +j_7 -j_8,
	j_1+j_6 -j_{10}}(2j_{10})
	 \nonumber\\
	&&\qquad \times
	{R}_{q,2j_1, 2j_2, j_1 +j_2 +j_6 -j_4; j_1 +j_6 -j_{10},
	 j_1 +j_5 -j_4}(2j_4)
	(e^{j_1}_{j_8- j_7} \ten e^{j_2}_{j_9 - j_8} \ten
	e^{j_3}_{j_4 - j_9})\nonumber \\
	&&=\\
	&&\sum_{j_8, j_9, j_{10}}
	{R}_{q,2j_1,2 j_2, j_1 +j_2 +j_7 -j_9;
	 j_1 +j_7 -j_8, j_1 +j_{10}-j_9}(2j_9)
	 {R}_{q,2j_1,2 j_3, j_1 +j_3 +j_{10} -j_4 ;j_1+j_{10}-j_9,
	j_1+j_5-j_4}(2j_4)
	\nonumber\\
	&&\qquad \times
	{R}_{q, 2j_2,2 j_3, j_2+j_3+j_7-j_5;j_2+j_7-j_{10},j_2+j_6-j_5}(2j_5)
	(e^{j_1}_{j_8- j_7} \ten e^{j_2}_{j_9 - j_8}\ten
	e^{j_3}_{j_4 - j_9}).\nonumber
\label{eq:YB2}
\end{eqnarray}
By use of the defining relation of $\wtilde{R}$ in \eqref{eq:ER5}, 
this turns out to be equivalent to the following identity:
\begin{eqnarray}
	&&\sum_{j_{10}}
	\wtilde{R}^{j_2 j_3}_{q,j_9 - j_8, j_4 - j_9 ; 
	j_4 - j_{10}, j_{10} - j_8}(j_4)
	 \wtilde{R}^{j_1 j_3}_{q,j_8 - j_7, j_{10} - j_8 ; 
	j_{10} - j_6, j_6 - j_7}(j_{10})
	\wtilde{R}^{j_1 j_2}_{q,j_{10} - j_6, j_4 - j_{10} ; 
	j_4 - j_5, j_5 - j_6}(j_4)	\nonumber \\
	&&=\\
	&&\sum_{j_{10}}
	\wtilde{R}^{j_1 j_2}_{q,j_8 - j_7, j_9 - j_8 ; 
	j_9 - j_{10}, j_{10} - j_7}(j_9)
	 \wtilde{R}^{j_1 j_3}_{q,j_9 - j_{10}, j_4 - j_9 ; 
	j_4 - j_5, j_5 - j_{10}}(j_4)
	\wtilde{R}^{j_2  j_3}_{q,j_{10} - j_7, j_5 - j_{10} ; 
	j_5 - j_6, j_6 - j_7}(j_5)
	\nonumber
\label{eq:YB3}
\end{eqnarray}
which yields, in virtue of \eqref{eq:ER5} the following known
 identity \cite[(6.19)]{KiR} satisfied by $q$-Racah coefficients: 
\begin{eqnarray}
	&&\sum_{j_{10}} (-1)^{-j_7 + j_8 + j_{10} + j_6}
	q^{(c_{j_7} - c_{j_8} - c_{j_{10}}
	- c_{j_6})/2} \nonumber\\
&&\qquad \qquad	{}_qW^{j_7,j_1,j_8}_{j_{10}-j_2,j_9-j_{10},j_9-j_2}(j_9)~
	{}_qW^{j_{10},j_1,j_9}_{j_5-j_3,j_4-j_5,j_4-j_3}(j_4)~
	{}_qW^{j_7,j_2,j_{10}}_{j_6-j_3,j_5-j_6,j_5-j_3} (j_5)\nonumber\\
&&	=\sum_{j_{10}} (-1)^{-j_4 + j_9 + j_{10} + j_5}
 	q^{(c_{j_4} - c_{j_{10}} - c_{j_9}
	- c_{j_5})/2} \nonumber\\
&&\qquad \qquad
	{}_qW^{j_8,j_2,j_9}_{j_{10}-j_3,j_4-j_{10},j_4-j_3}(j_4)~
	{}_qW^{j_7,j_1,j_8}_{j_6-j_3,j_{10}-j_6,j_{10}-j_3}(j_{10})~
	{}_qW^{j_6,j_1,j_{10}}_{j_5-j_2,j_4-j_5,j_4-j_2}(j_4).
\label{eq:YB4}
\end{eqnarray}
Then, by use of \eqref{eq:YB4} with \eqref{eq:ER9} substituted
and symmetries of $q$-$6j$ symbols (i.e., the $q$-$6j$ symbol is invariant
under any permutation of 
columns and also under an interchange of upper and lower arguments in each
of any two of its columns, see \cite{N1989}), we finally show that 
$q$-$6j$ symbols satisfy the following identity:
\begin{eqnarray}
 	&&\sum_{j_{10}} (-1)^{2j_{10}} [2j_{10} + 1]_{q}~ 
	q^{-(c_{j_{10}}+c_{j_8}+c_{j_6}+c_{j_4})/2}\nonumber\\
	&&\qquad\qquad\times
\qsixj{j_2}{j_7}{j_{10}}{j_1}{j_9}{j_8}
\qsixj{j_2}{j_6}{j_5}{j_3}{j_{10}}{j_7}
\qsixj{j_1}{j_5}{j_4}{j_3}{j_9}{j_{10}}
\nonumber\\
	&&=\sum_{j_{10}} (-1)^{2j_{10}} [2j_{10} + 1]_{q}~ 
	q^{-(c_{j_{10}}+c_{j_7}+c_{j_9}+c_{j_5})/2} 
	\nonumber\\
	&&\qquad\qquad\times
\qsixj{j_1}{j_6}{j_{10}}{j_3}{j_8}{j_7}
\qsixj{j_2}{j_{10}}{j_4}{j_3}{j_9}{j_8}
\qsixj{j_2}{j_6}{j_5}{j_1}{j_4}{j_{10}}.
\label{eq:YB5}
\end{eqnarray}

\quad\\
\begin{footnotesize}
\begin{quote}
{ T. H. Koornwinder, Korteweg-de Vries Institute, University of
 Amsterdam,\\
 Plantage Muidergracht 24, 1018 TV Amsterdam, The Netherlands;

\vspace{\smallskipamount}
email: }{\tt thk@science.uva.nl}

\vspace{\bigskipamount}
{ N. Touhami, Korteweg-de Vries Institute, University of
Amsterdam,\\
Plantage Muidergracht 24, 1018 TV Amsterdam, The Netherlands; 

\vspace{\smallskipamount}
Laboratoire de Physique Th\'eorique, Universit\'e d'Oran Es-S\'enia, 
Oran, Algeria;

\vspace{\smallskipamount}
email: }{\tt touhami@science.uva.nl}
\end{quote}
\end{footnotesize}
\end{document}